\documentclass[11pt]{article}
\usepackage{auto-pst-pdf}
\usepackage{latexsym,color,graphicx,amsmath,amssymb,amsthm}
\usepackage{url}

\newtheorem{theorem}{Theorem}[section]
\newtheorem{lemma}[theorem]{Lemma}
\newtheorem{proposition}[theorem]{Proposition}
\newtheorem{corollary}[theorem]{Corollary}

\def\reals{{\mathbb R}}
\def\cplx{{\mathbb C}}
\def\P{{\mathbb P}}

\def\sph{{\mathbb S}}

\def\C{{\mathcal C}}

\begin{document}

\def\reals{{\mathbb R}}
\def\cplx{{\mathbb C}}
\def\deg{{\mathsf{deg}}}

\newcommand{\pts}{\mathcal P}
\newcommand{\vrts}{\mathcal V}
\newcommand{\curves}{\mathcal C}
\newcommand{\lines}{\mathcal L}
\newcommand{\motions}{\mathcal M}
\newcommand{\planes}{\Pi}
\newcommand{\flats}{\mathcal F}
\newcommand{\DS}{\ensuremath{\mathbb S}}
\newcommand{\dist}[1]    {{\,{\mathrm{D}_3}\!\left(#1\right)}}
\newcommand{\res}{\text{Res}}

\newcommand{\Deg}{D} 
\makeatletter
\newcommand{\fl}[1]{\mathsf{FL}_{#1}}
\newcommand{\pfl}[1]{\mathsf{PFL}_{#1}}
\newcommand{\ProofEndBox}{{\ifhmode\unskip\nobreak\hfil\penalty50 \else
          \leavevmode\fi\quad\vadjust{}\nobreak\hfill$\Box$
            \finalhyphendemerits=0 \par}}
\makeatother
\newcommand{\proofend}{\ProofEndBox\smallskip}
\title{Incidences between points and lines on two- and three-dimensional varieties\footnote{%
Work on this paper by Noam Solomon and Micha Sharir was supported by
Grant 892/13 from the Israel Science Foundation. Work by Micha
Sharir was also supported by Grant 2012/229 from the U.S.--Israel
Binational Science Foundation, by the Israeli Centers of Research
Excellence (I-CORE) program (Center No.~4/11), by the Blavatnik Computer Science Research Fund at Tel Aviv University,
and by the Hermann Minkowski-MINERVA Center for Geometry at Tel Aviv University. An earlier version of the paper, whoch only contains some of
the results and in a weaker form, is:
M. Sharir and N. Solomon,
Incidences between points and lines on a two-dimensional variety,
in arXiv:1501.01670. }}

\author{
Micha Sharir\thanks{%
School of Computer Science, Tel Aviv University, Tel Aviv 69978,
Israel. {\sl michas@post.tau.ac.il} } \and
Noam Solomon\thanks{%
School of Computer Science, Tel Aviv University, Tel Aviv 69978,
Israel. {\sl noam.solom@gmail.com} } }

\maketitle

\begin{abstract}
Let $P$ be a set of $m$ points and $L$ a set of $n$ lines in $\reals^4$, 
such that the points of $P$ lie on an algebraic three-dimensional surface 
of degree $\Deg$ that does not contain hyperplane or quadric components, 
and no 2-flat contains more than $s$ lines of $L$. 
We show that the number of incidences between $P$ and $L$ is
$$
I(P,L) = O\left(m^{1/2}n^{1/2}D + m^{2/3}n^{1/3}s^{1/3} + nD + m\right) ,
$$
for some absolute constant of proportionality.
This significantly improves the bound of the authors~\cite{SS4d}, for arbitrary
sets of points and lines in $\reals^4$, when $\Deg$ is
not too large. Moreover, when $\Deg$ and $s$ are constant, we get a
\emph{linear} bound. The same bound holds when the three-dimensional
surface is embedded in any higher dimensional space.

\smallskip

For the proof of this bound, we revisit certain parts of~\cite{SS4d}, combined 
with the following new incidence bound, for which we present a direct and fairly 
simple proof. Let $P$ be a set of $m$ points and $L$ a set of $n$ lines in $\reals^d$, 
for $d\ge 3$, which lie in a common two-dimensional algebraic surface of degree $\Deg$
(assumed to be $\ll n^{1/2}$) that does not contain any 2-flat, so that no 2-flat contains 
more than $s$ lines of $L$ (here we require that the lines of $L$ also be contained in the surface). 
Then the number of incidences between $P$ and $L$ is
$$
I(P,L) = O\left(m^{1/2}n^{1/2}\Deg^{1/2} + m^{2/3}\Deg^{2/3}s^{1/3} + m + n\right).
$$
When $d=3$, this improves the bound of Guth and Katz~\cite{GK2} for this special
case, when $\Deg\ll n^{1/2}$. Moreover, the bound does not involve the term $O(nD)$.
This term arises in most standard approaches, and its removal is a significant aspect of our
result. Again, the bound is linear when $\Deg = O(1)$.

\smallskip

Finally, we also obtain (slightly weaker) variants of both results over the complex field.
For two-dimensional varieties, when the ambient space is $\cplx^3$
(or any higher-dimensional $\cplx^d$), under the same assumptions as above, we have
$$
I(P,L) = O\left(m^{1/2}n^{1/2}\Deg^{1/2} + m^{2/3}\Deg^{2/3}s^{1/3} + \Deg^3 + m + n\right).
$$
For three-dimensional varieties, embedded in $\cplx^4$
(or any higher-dimensional $\cplx^d$), under the same assumptions as above, we have
$$
I(P,L) = O\left(m^{1/2}n^{1/2}D + m^{2/3}n^{1/3}s^{1/3} + D^6 + m + nD \right) .
$$
These new incidence bounds are among the very few bounds, known so far, that hold over the complex field.
The bound for two-dimensional (resp., three-dimensional) varieties coincides with the bound in the real case when 
$\Deg = O(m^{1/3})$ (resp., $\Deg = O(m^{1/6})$).
\end {abstract}

\section{Introduction}

Let $P$ be a set of $m$ distinct points and $L$ a set of $n$
distinct lines in $\reals^d$ or in $\cplx^d$. Let $I(P,L)$ denote the number of
incidences between the points of $P$ and the lines of $L$; that is,
the number of pairs $(p,\ell)$ with $p\in P$, $\ell\in L$, and
$p\in\ell$. If all the points of $P$ and all the lines of $L$ lie in
a common 2-flat, then, in the real case, the classical Szemer\'edi--Trotter
theorem~\cite{SzT} yields the worst-case tight bound
\begin{equation} \label{inc2}
I(P,L) = O\left(m^{2/3}n^{2/3} + m + n \right) .
\end{equation}
The same bound also holds in the complex plane, as shown later by T\'oth~\cite{Tot} and Zahl~\cite{Zah}.

This bound clearly also holds in $\reals^d$ and in $\cplx^d$, for any $d$, by projecting the given
lines and points onto some generic 2-flat. Moreover, the bound will
continue to be worst-case tight by placing all the points and lines
in a common 2-flat, in a configuration that yields the planar lower
bound.

In the 2010 groundbreaking paper of Guth and Katz~\cite{GK2}, an improved
bound for the real case has been derived for $I(P,L)$, for a set $P$ of $m$ points and a set
$L$ of $n$ lines in $\reals^3$, provided that not too many lines of $L$ lie
in a common plane. Specifically, they showed:\footnote{%
  We skip over certain subtleties in their bound: They also assume that no
  \emph{regulus} contains more than $s=\sqrt{n}$ input lines, but then they are able
  also to bound the number of intersection points of the lines. Moreover,
  if one also assumes that each point is incident to at least three lines
  then the term $m$ in the bound can be dropped.}

\begin{theorem}[Guth and Katz~\cite{GK2}]
\label {ttt}
Let $P$ be a set of $m$ distinct points and $L$ a set of $n$ distinct lines
in $\reals^3$, and let $s\le n$ be a parameter,
such that
no plane contains more than $s$ lines of $L$. Then
$$
I(P,L) = O\left(m^{1/2}n^{3/4} + m^{2/3}n^{1/3}s^{1/3} + m + n\right).
$$
\end{theorem}
This bound (or, rather, an alternative formulation thereof) was a
major step in the derivation of the main result of \cite{GK2}, which
was an almost-linear lower bound on the number of distinct
distances determined by any finite set of points in the plane, a
classical problem posed by Erd{\H o}s in 1946~\cite{Er46}. Guth and Katz's
proof uses several nontrivial tools from algebraic and differential
geometry, most notably the Cayley--Salmon--Monge theorem on osculating
lines to algebraic surfaces in $\reals^3$, and various properties
of ruled surfaces. All this machinery comes on top of the major innovation 
of Guth and Katz, the introduction of the \emph{polynomial partitioning technique}. 

For the purpose of the analysis in this paper, it is important to recall, 
right away, that the polynomial partitioning technique holds only over the reals. 
This will be the major stumbling block that we will face as we handle the complex case.
We overcome (or rather bypass) it by exploiting the fact that all the lines are also 
contained in the given variety; see below for details.

In four dimensions, and for the real case, the authors established in \cite{SS4d} 
a sharper and (almost) optimal bound. More precisely, they have shown:
\begin{theorem} \label{th:mainfocs}
Let $P$ be a set of $m$ distinct points and $L$ a set of $n$
distinct lines in $\reals^4$, and let $s\le q\le n$ be
parameters, such that (i) no hyperplane or quadric contains more than
$q$ lines of $L$, and (ii) no 2-flat contains more than $s$ lines of $L$. Then,
\begin {equation}
\label{ma:in0}
I(P,L) \le 2^{c\sqrt{\log m}} \left( m^{2/5}n^{4/5} +
m \right) + A\left( m^{1/2}n^{1/2}q^{1/4} + m^{2/3}n^{1/3}s^{1/3} +
n\right) ,
\end {equation}
for suitable (absolute) constant parameters $A$ and $c$.
Moreover, except for the factor $2^{c\sqrt{\log m}}$, the bound is
tight in the worst case, for all $m$ and $n$ and suitable ranges of $q$ and $s$.
For certain ranges of $m$ and $n$ the bound holds without that factor.
\end{theorem}

\paragraph{Our results.}
In this work we ``contest'' the leading terms
$O\left(m^{1/2}n^{3/4}\right)$ (for $d=3$) and
$2^{c\sqrt{\log m}} m^{2/5}n^{4/5}$ (for $d=4$),
and show that in many situations they can be significantly improved.
See item (2) in the discussion in Section~\ref{Se:dis} for further elaboration of this issue.
A major feature of this work is that, in the setups considered here, the analysis
can also be carried over to the complex domain, except for a small penalty that we 
pay for bypassing the polynomial partitioning technique, which, as noted, only holds over the reals.

Concretely, we assume that the points of $P$ lie on some 
algebraic variety, and derive significantly improved bounds when the degree of the variety is not too large.
In the former case we assume that the points and the lines lie on a two-dimensional variety,
which is allowed to be embedded in any $\reals^d$, for $d\ge 3$.
In the latter case we assume that the points (but not necessarily the lines) lie on a three-dimensional
variety, embedded in any $\reals^d$, for $d\ge 4$. In the former (resp., latter) 
case we also assume that the variety contains no plane (resp., no hyerplane or quadric).
Thus, in addition to improving the respective bounds in Theorems~\ref{ttt}
and \ref{th:mainfocs}, for the special cases under consideration, and in extending them to the complex domain,
we obtain an extra bonus by extending the results to two-dimensional and
three-dimensional varieties embedded in any higher dimension.

\paragraph{Points on a two-dimensional variety.} 
We derive two closely related results, one that holds over the real field and one
that holds also over the complex field.
It is simplest to think of the variety as embedded in $\reals^3$ or in $\cplx^3$.
The real case is a special case of the setup of Guth and Katz~\cite{GK2}, where there is no
need to use the polynomial partitioning method, because we assume
that the points and lines all lie in a common surface (the zero set of a polynomial)
of degree $\Deg$. This very assumption is also the one that lets us derive the (slightly weaker) version
that holds over $\cplx$, thereby constituting a significant progress over the existing 
theory of incidences in three (and higher) dimensions. (To be more precise, over the reals
we do apply the polynomial partitioning technique (as a step in the application of the Guth-Katz bound), 
but only to a small subset of the lines.)

Concretely, our first main result, for this setup, is the following theorem.

\begin{theorem} \label {th:main3d} 
{\bf (a) The real case:}
Let $P$ be a set of $m$ distinct points and $L$ a
set of $n$ distinct lines in $\reals^d$, for any $d\ge 3$, and let $2\le s \le \Deg$
be two integer parameters, so that all the points and
lines lie in a common two-dimensional algebraic variety $V$ of
degree $\Deg$ that does not contain any 2-flat, and so that no
2-flat contains more than $s$ lines of $L$. Then
\begin{equation} \label{st}
I(P,L) = O\left(m^{1/2}n^{1/2}\Deg^{1/2} + m^{2/3}\Deg^{2/3}s^{1/3} + m + n\right) .
\end{equation}

\smallskip
\noindent
{\bf (b) The complex case:}
Under exactly the same assumptions, when the ambient space is $\cplx^d$, for any $d\ge 3$,
we have
\begin{equation} \label{stc}
I(P,L) = O\left(m^{1/2}n^{1/2}\Deg^{1/2} + m^{2/3}\Deg^{2/3}s^{1/3} + \Deg^3 + m + n\right) .
\end{equation}
\end{theorem}

The assumption that $s$ is at most $\Deg$ can be dropped, because, for any
2-flat $\pi$, the intersection $\pi \cap V$ is a one-dimensional
plane algebraic curve of degree at most $\Deg$ in $\pi$ (this holds since
$V$ does not contain any 2-flat), and can therefore contain at most $\Deg$ lines.

We also have the following easy and interesting corollary.
\begin{corollary}
\label {cor:main} Let $P$ be a set of $m$ distinct points and $L$ a
set of $n$ distinct lines in $\reals^d$ or in $\cplx^d$, for any $d\ge 3$, such that all the points and
lines lie in a common two-dimensional algebraic variety of constant
degree that does not contain any 2-flat. Then
$I(P,L) = O\left(m + n\right)$,
where the constant of proportionality depends on the degree of the surface.
\end{corollary}

For $d=3$, the corollary can also be derived, for the real case, from
the analysis in Guth and Katz~\cite{GK2}, using a somewhat different
approach. Moreover, although not explicitly stated there, it seems that
the argument in \cite{GK2} also works over $\cplx$.
As a matter of fact, the corollary can also be extended (with a different bound though)
to the case where the containing surface may have planar components. 
See a remark to that effect in Corollary~\ref{cor:mainx} in the concluding section.

We also exploit the proof technique of Theorem~\ref{th:main3d} to derive an upper bound
of $O(n\Deg)$ on the number of \emph{$2$-rich points} determined by a set of $n$ lines
contained in a variety, as above, in both the real and complex cases. 
See Section~\ref{Se:dis} for details.

We note that if we only assume that the points lie on $V$ but the lines are arbitrary, we will
incur the term $O(nD)$, which we are trying to avoid. In the case of three-dimensional varieties, 
discussed below, we cannot avoid this term even when the lines lie on the variety.

The significance of Theorem~\ref{th:main3d} is fourfold:

\smallskip

\noindent (a) First and foremost, the theorem yields a new incidence result for points and lines on a 
two-dimensional variety \emph{over the complex field}, in three and higher dimensions. 
Incidence results over the complex domain are rather rare. They include (as already mentioned)
T\'oth's extension of the Szemer\'edi-Trotter bound to the complex plane~\cite{Tot},
which was the only result of that kind that predated the introduction of the algebraic machinery
by Guth and Katz, and several more recent works~\cite{ShZa15,SoTa,SZ15,Zah} 
(where the latter work~\cite{Zah} provides an alternative algebraic derivarion of T\'oth's bound).

\smallskip

\noindent (b) In the real three-dimensional case, the bound improves the Guth--Katz
bound when $D\ll n^{1/2}$, for two-dimensional varieties $V$ that do
not contain planes. Note that the threshold $n^{1/2}$ is a natural
one because, as is well known and easy to show, any set of $n$ lines
in $\reals^3$ admit a polynomial of degree $O(n^{1/2})$ whose zero
set contains all the lines; a simple modification of the construction
applies in higher dimensions too. (Of course the comparison is far from perfect,
because this polynomial may have many linear components, in which case our bound does not apply.
Still, it offers some basis for evaluating the quality of our bound.)
In three dimensions, this threshold is in fact larger than the standard degree
$O(m^{1/2}/n^{1/4})$ used in the analysis of Guth and
Katz~\cite{GK2}, when $m<n^{3/2}$.

\smallskip

\noindent
(c) Another significant feature of our bound is that it does not contain
the term $nD$, which arises naturally in \cite{GK2} and other works,
and seems to be unavoidable when $P$ is an arbitrary set of points. 
When $D$ is not a constant, this becomes a crucial feature of the new bound, which
has already been exploited in the analysis in~\cite{SS4d},
and is also used in the second main result of this paper, Theorem~\ref{th:main4d} below.\footnote{%
  Although the bound in Theorem~\ref{th:main4d} does contain the term $nD$, it still 
  crucially relies on the absence of this term in the bound for two-dimensional varieties.}
See an additional discussion of this feature at the end of the paper.

\smallskip

\noindent
(d) It offers a sharper point-line incidence bound in arbitrary dimensions,
for the special case assumed in the theorem (which again holds over the complex field too).

Theorem~\ref{th:main3d}(a) has been used, as one of the key tools, in the
analysis in our paper~\cite{SS4d} on incidences between points and
lines in four dimensions (see below). In this application, the absence of the
term $nD$ is a crucial feature of our result, which was required in the scenario considered in \cite{SS4d}.

The proof of Theorem~\ref{th:main3d} makes extensive use of several properties of ruled surfaces
in $\reals^3$ or in $\cplx^3$. While these results exist as folklore in the literature,
and short proofs are provided for some of them, e..g., in \cite{GK2},
we include here detailed and rigorous proofs thereof, making them more accessible
to the combinatorial geometry community. Other recent expositions include 
Guth's recent survey~\cite{Gut:surv} and book~\cite{Gut:book}, 
and a survey by Koll\'ar~\cite{Kollar}.\footnote{%
  While there is (naturally) some overlap between these surveys and our exposition, the main
  technical properties that we present do not seem to be rigorously covered in the other works.}

\paragraph{Points on a three-dimensional variety.}
Our second main result deals with the case where the points lie on a
three-dimensional variety, embedded in $\reals^4$ or in $\cplx^4$, or in any higher dimension.
Similar to the three-dimensional case discussed above, we have to be careful here too, 
because hyperplanes and 3-quadrics (in $\reals^4$, and, aposteriorily, in $\cplx^4$ too) 
admit ``too many'' incidences in the worst case. That is, by a generalization of Elekes's
construction~\cite{Elekes}, there exists a configuration of $m$
points and $n$ lines on a 3-flat with $\Theta(m^{1/2}n^{3/4})$
incidences. More recently, Solomon and Zhang~\cite{SZ} established an
analogous statement for three-dimensional quadrics, when
$n^{9/8}<m<n^{3/2}$. Concretely, for such values of $m$ and $n$,
they have constructed a quadric $S\subset\reals^4$, a set $P$ of $m$
points on $S$, and a set $L$ of $n$ lines contained in $S$, so that
(i) the number of lines in any common 2-flat is $O(1)$, (ii) the
number of lines in any hyperplane is $O(n/m^{1/3})$, and (iii) the
number of incidences between the points and lines is
$\Omega(m^{2/3}n^{1/2})$, which is asymptotically larger than the
corresponding bound in (\ref{ma:in0}) (for $s=O(1)$), when
$n^{9/8} \ll m\le n^{3/2}$.

In other words, when studying incidences with points on a variety in $\reals^4$ or in $\cplx^4$,
the cases where the variety is (or contains) a hyperplane or a quadric are special and do not yield
sharper bounds. (In the real case, the case of a hyperplane puts us back in $\reals^3$,
where the best bound is Guth and Katz's in Theorem~\ref{ttt}, and the case of a quadric reduces to
the same setup via a suitable generic projection onto $\reals^3$.)
In the second main result of this paper, we show that if all the points
lie on a three-dimensional algebraic variety of degree $D$ without
3-flat or 3-quadric components, and if no 2-flat contains more than
$s$ lines, then, if $D$ and $s$ are not too large, the bound becomes
significantly smaller. Moreover, here too we get a real version and a complex version
of the theorem. Specifically, we show:

\begin{theorem} \label{th:main4d}
{\bf (a) The real case:}
Let $P$ be a set of $m$ distinct points and $L$ a set of $n$
distinct lines in $\reals^d$, for any $d \ge 4$, and let $s$ and $D$
be parameters, such that
(i) all the points of $P$ lie on a three-dimensional algebraic
variety of degree $\Deg$, without 3-flat or 3-quadric components,
and (ii) no 2-flat contains more than $s$ lines of $L$. Then
\begin {equation}
\label {ma:in} I(P,L) = O\left( m^{1/2}n^{1/2}D + m^{2/3}n^{1/3}s^{1/3} + nD + m\right) .
\end {equation}
When $D$ and $s$ are constants, we get the linear bound $O(m+n)$.

\smallskip
\noindent
{\bf (b) The complex case:}
Under exactly the same assumptions, when the ambient space is $\cplx^d$, for any $d\ge 4$,
we have
\begin{equation} 
\label {ma:incplx} I(P,L) = O\left( m^{1/2}n^{1/2}D + m^{2/3}n^{1/3}s^{1/3} + \Deg^6 + nD + m\right) .
\end{equation}
\end{theorem}

\medskip

\noindent {\bf Remarks. (1)} 
Note that for $D < \min\{m^{1/2}/n^{1/4},\;n^{1/4}\}$, our bound for the real case is sharper than the bound of 
Guth and Katz~\cite{GK2} (note that $m^{1/2}/n^{1/4}$ is the degree of the partitioning 
polynomial used in the analysis of \cite{GK2} for $m\le n^{3/2}$). On the other hand, when 
$D > \min\{m^{1/2}/n^{1/4},\;n^{1/4}\}$, our bound is not the best possible. Indeed, in this case we can project
$P$ and $L$ onto some generic $3$-flat, and apply instead the bound of~\cite{GK2} to the
projected points and lines (we also show that a generic choice of the image $3$-flat also ensures that no 
$2$-flat contains more than $s$ of the projected lines), which is sharper than ours for these values of $D$. 
In the complex case, we also need to assume that $D$ is small enough so that the term $D^6$ does not dominate the other terms,
so as to make the bound look like the bound in the real case.
\ref{ma:incplx} 

\medskip

\noindent {\bf (2)} 
As already noted earlier, here we do not have to insist that the lines of $L$ be contained in the variety.
A line not contained in a variety of degree $D$ can intersect it at most $D$ time, so the number of
incidences with such lines is at most $nD$. (The actual argument that yields this term is more involved, because we apply
the argument to a variety of larger degree; see below for details.)
For two-dimensional varieties, since we want to avoid this term,
we require the lines to be contained in the variety. Here, since we allow this term, the lines can be arbitrary.

\medskip

\noindent {\bf (3)} 
The assumption that the points of $P$ lie on a variety is not as restrictive as it might sound, because,
in four dimensions,
one can always construct a polynomial $f$ of degree $O(m^{1/4})$ whose zero set contains all the given points,
or alternatively, a polynomial $f$ of degree $O(n^{1/3})$ whose zero set fully contains all the given lines.
The assumptions become restrictive, and the bound becomes more interesting, when $D$ is significantly smaller.
In addition, the constructed polynomial $f$ could be one for which $Z(f)$ does contain 3-flats or 3-quadrics,
so another restrictive aspect of our assumptions is that they exclude these situations.

Similar to the case of a two-dimensional variety, we also have the following easy and interesting corollary
(it does not hold when $V$ contains 3-flats or 3-quadrics, in light of the lower bound constructions in 
\cite{SS4d,SZ}).
\begin{corollary} \label{cor:main4d} 
Let $P$ be a set of $m$ distinct points and $L$ a set of $n$ distinct lines in 
$\reals^d$ or in $\cplx^d$, for any $d\ge 4$, such that all the points
and lines lie in a common three-dimensional algebraic variety $V$ of
constant degree that does not contain any 3-flats or 3-quadrics,
and no 2-flat contains more than $O(1)$ lines of $L$. Then
$I(P,L) = O\left(m + n\right)$, where the constant of
proportionality depends on the degree of $V$.
\end{corollary}

Theorem~\ref{th:main3d} is a key technical ingredient in the proof
of Theorem~\ref{th:main4d}. The proofs of both theorems are somewhat
technical, and use a battery of sophisticated tools from algebraic
geometry. Some of these tools are borrowed and adapted from our
previous work~\cite{SS4d}. Other tools involve properties of ruled
surfaces, which, as already said, are established here rigorously, for the sake of completeness.
Since most of the presentation and derivation of these results is
within the scope of algebraic and differential geometry, we delegate 
most of this part to an appendix, keeping only a brief overview in the following section.

The proof of Theorem~\ref{th:main3d} is then presented in
Section~\ref{ssec:pfa}, and the proof of Theorem~\ref{th:main4d} is
presented in Section \ref{sec:pfmain}. The concluding Section~\ref{Se:dis}
discusses our results, establishes a few consequences thereof, and raises several related open problems.

\section{Algebraic preliminaries and ruled surfaces} \label{Se:pre}

In this section we present a brief overview of ruled surfaces and some of the algebraic machinery. 
Full details are given in the appendix.

\paragraph{Ruled surfaces.} 
We say that a real (resp., complex) surface $V$ is \emph{ruled by real} 
(resp., \emph{complex}) \emph{lines} if every point $p\in V$ is incident 
to a real (resp., complex) line that is fully contained in $V$. 
The surface is \emph{doubly ruled} if each point is incident to two distinct lines that are fully contained in $V$.
See the appendix for a slightly looser definition and more details, including the not-so-extensive
literature on ruled surfaces. The recent survey~\cite{Gut:surv} and book~\cite{Gut:book},
as well as an exposition in~\cite{Kollar} make up for this sparsity.

\paragraph{Flecnodes in three dimensions and the Cayley--Salmon--Monge Theorem.}
We first recall the classical theorem of Cayley and Salmon, also due to Monge. Consider
a polynomial $f \in \cplx[x,y,z]$ of degree $D\ge 3$. A \emph{flecnode} 
of $f$ is a point $p$ on the zero set $Z(f)$ of $f$, for which there exists
a line that is incident to $p$ and \emph{osculates} to $Z(f)$ at $p$
to order three. That is, if the direction of the line is $v$ then
$f(p) = 0$, and $\nabla_v f(p) = \nabla_v^2 f(p) = \nabla_v^3 f(p) = 0$, 
where $\nabla_v f, \nabla^2_v f, \nabla^3_v f$ are, respectively, the 
first, second, and third-order derivatives of $f$ in the direction $v$.
The \emph{flecnode polynomial} of $f$, denoted $\fl{f}$, is the polynomial 
obtained by eliminating $v$, via resultants, from these three homogeneous equations 
(where $p$ is regarded as a fixed parameter). As shown in 
Salmon~\cite[Chapter XVII, Section III]{salmon}, the degree
of $\fl{f}$ is at most $11D-24$. By construction, the flecnode
polynomial of $f$ vanishes on all the flecnodes of $f$, and in
particular on all the lines fully contained in $Z(f)$.

\begin{theorem}[Cayley and Salmon~\cite{salmon}, Monge~\cite{Mon}]
\label{th:flec2a} Let $f \in \cplx[x,y,z]$ be a polynomial of degree
$\Deg\ge 3$. Then $Z(f)$ is ruled by (complex) lines if and only if
$Z(f) \subseteq Z(\fl{f})$.
\end{theorem}

The notions of flecnodes and of the flecnode polynomial can be extended to
four dimensions, as done in~\cite{SS4d}. Informally, the 
\emph{four-dimensional flecnode polynomial} $\fl{f}^4$ of a 4-variate polynomial 
$f$ is defined analogously to the three-dimensional variant $\fl{f}$, 
and captures the property that a point on $Z(f)$ is incident to 
a line that osculates to $Z(f)$ up to the \emph{fourth} order.
It is obtained by eliminating the direction $v$ of the osculating
line, using resultants, from the four homogeneous equations given by the vanishing of
the first four terms of the Taylor expansion of $f(p+tv)$ near $p$.
Clearly, $\fl{f}^4$ vanishes identically on every line that is fully contained
in $Z(f)$. As in the three-dimensional case, its degree can be shown to be $O(D)$.

Landsberg~\cite{Land} derives an analog of Theorem~\ref{th:flec2a} that holds 
for three-dimensional surfaces (see~\cite[Theorem 2.11]{SS4d}). Specifically, 
Landsberg's theorem asserts that if $\fl{f}^4$ vanishes identically on $Z(f)$, 
then $Z(f)$ is ruled by (possibly complex) lines. We will discuss this
in more detail in Section~\ref{sec:pfmain}. These theorems, in three and 
four dimensions, play an important role in the proofs of the main theorems.


\paragraph{Reguli.} 
A \emph{regulus} is the surface (in $\reals^3$ or $\cplx^3$) spanned by all lines 
that meet three pairwise skew lines in $3$-space.\footnote{%
  Technically, in some definitions (cf., e.g., Edge~\cite[Section I.22]{Edge})
  a regulus is a one-dimensional family of generator lines of the actual
  surface, i.e., a curve in the Pl\"ucker or Grassmannian space of lines,
  but we use here the alternative notion of the surface spanned by these lines.}
For an elementary proof that a doubly ruled surface over $\reals$ must be a
regulus, we refer the reader to Fuchs and Tabachnikov~\cite[Theorem 16.4]{FT}.
We give in the appendix a proof of the following more general characterization of 
doubly ruled surfaces in $\reals^3$ or $\cplx^3$. 
\begin {lemma}
\label{doubly} Let $V$ be an irreducible ruled surface in $\reals^3$
or in $\cplx^3$ which is not a plane, and let $\C\subset V$ be an
algebraic curve, such that every non-singular point $p\in V\setminus \C$ 
is incident to exactly two lines that are fully contained in $V$. Then $V$ is a regulus.
\end {lemma}

\paragraph{Singly ruled surfaces.}
Ruled surfaces that are neither planes nor reguli are called
\emph{singly ruled} surfaces (a terminology justified by
Theorem~\ref{singly}, given below). A line $\ell$, fully contained
in an irreducible singly ruled surface $V$, such that every point of
$\ell$ is ``doubly ruled'', i.e., every point on $\ell$ is incident
to another line fully contained in $V$, is called an
\emph{exceptional} line of $V$. A point $p_V \in V$ that is incident 
to infinitely many lines fully contained in
$V$ is called an \emph{exceptional} point of $V$.

The following result is
another folklore result in the theory of ruled surfaces, used in
many studies (such as Guth and Katz~\cite{GK2}).
It justifies the terminology ``singly-ruled surface'', by showing that
the surface is generated by a one-dimensional family of lines, and that each point
on the surface, with the possible exception of points lying on some curve,
is incident to exactly one generator. It also shows that there are only finitely
many exceptional lines; the property that their number is at most two (see~\cite{GK2})
is presented later. We give (in the appendix) a detailed and rigorous proof, to make our presentation
as self-contained as possible; we are not aware of any similarly detailed argument in the literature.
\begin{theorem}
\label{singly} (a) Let $V$ be an irreducible ruled two-dimensional
surface of degree $\Deg>1$ in $\reals^3$ (or in $\cplx^3$), which is
not a regulus. Then, except for finitely many exceptional lines, the
lines that are fully contained in $V$ are parameterized by an
irreducible algebraic curve $\Sigma_0$ (in the parametric Pl\"ucker 
space $\P^5$ that represents lines in 3-space; see the appendix), 
and thus yield a 1-parameter family of generator lines $\ell(t)$, for 
$t\in \Sigma_0$, that depend continuously on the real or complex parameter 
$t$. Moreover, if $t_1 \ne t_2$, and $\ell(t_1) \ne \ell(t_2)$, then there 
exist sufficiently small and disjoint neighborhoods $\Delta_1$ of $t_1$ and 
$\Delta_2$ of $t_2$, such that all the lines $\ell(t)$, for 
$t\in \Delta_1\cup \Delta_2$, are distinct.

\smallskip

\noindent (b) There exists a one-dimensional curve $\C\subset V$, such that 
any point $p$ in $V\setminus\C$ is incident to exactly one line fully contained in $V$.
\end{theorem}

\paragraph{Exceptional lines on a singly ruled surface.}
In view of Theorem~\ref{singly},
every point on a singly ruled surface $V$ is incident to at least
one generator. Hence an exceptional (non-generator) line is a line
$\ell\subset V$ such that every point on $\ell$ is incident to a
generator (which is different from $\ell$). 
We establish in the appendix the following property.
\begin {lemma} \label{only2exc}
Let $V$ be an irreducible ruled surface in $\reals^3$ or in $\cplx^3$, 
which is neither a plane nor a regulus. Then (i) $V$ contains at most 
two exceptional lines, and (ii) $V$ contains at most one exceptional point.
\end {lemma}

\paragraph{Generic projections preserve non-planarity.}
In the analysis in Section~\ref{ssec:pfa}, the goal is to project
$\reals^d$ onto some generic 3-flat so that non-coplanar triples of
lines do not project to coplanar triples. This is easily achieved by
repeated applications of the following technical result, reducing
the dimension one step at a time.
\begin {lemma}
\label{ap:th} Let $\ell_1,\ell_2,\ell_3$ be three non-coplanar lines
in $\reals^d$, for $d\ge 4$. Then, under a generic projection of $\reals^d$ onto
some hyperplane $H$, the respective images
$\ell^*_1,\ell^*_2,\ell^*_3$ of these lines are still non-coplanar.
\end {lemma}

\paragraph{Real vs.~complex.}
With a few exceptions, our analysis holds over the complex field too, and therefore
most of the unfolding analysis is carried out over $\cplx$, the main reason being that 
most of the basic tools in algebraic geometry only hold over $\cplx$ 
(or other algebraically closed fields).
There are a few places, though, where the analysis holds only over the reals, such as
in the (explicit or implicit) applications of the polynomial partitioning technique
of Guth and Katz~\cite{GK2}. Passing to the complex domain (and sometimes also to the
projective setting) does not pose any difficulties for upper bounds in real incidence 
problems---every real incidence will be preserved, and at worst we will be counting 
additional incidences, on the non-real portion of the extended varieties. 
With this understanding, and with the appropriate caution, we will move freely between
the real and complex domains, as convenient, staying mostly in the complex domain.

\section{Proof of Theorem~\ref{th:main3d}} \label{ssec:pfa}

In most of the analysis in this section, we will consider the case
$d=3$. The reduction from an arbitrary dimension to $d=3$ will be
presented at the end of the section.

We will prove both parts of the theorem ``hand in hand'', bifurcating 
(in a significant manner) only towards the end of the analysis.

For a point $p$ on an irreducible singly ruled surface $V$, which is not the 
exceptional point of $V$ (see Section~\ref{Se:pre} and the appendix for its 
definition), we let $\Lambda_V(p)$ denote the number of generator lines 
incident to $p$ and fully contained in $V$ (so if $p$ is incident to an 
exceptional line, we do not count that line in $\Lambda_V(p)$). We also put 
$\Lambda_V^*(p) := \max\{0,\Lambda_V(p)-1\}$. Finally, if $V$ is a cone and 
$p_V$ is its exceptional point (that is, apex), we put 
$\Lambda_V(p_V) = \Lambda_V^*(p_V):=0$. We also consider a variant of this notation,
where we are also given a finite set $L$ of lines (where not all
lines of $L$ are necessarily contained in $V$), which does not
contain any of the (at most two) exceptional lines of $V$.  For a
point $p\in V$, we let $\lambda_V(p;L)$ denote the number of lines
in $L$ that pass through $p$ and are fully contained in $V$, with
the same provisions as above, namely that we do not count incidences
with exceptional lines, nor do we count incidences occurring at an
exceptional point, and put $\lambda_V^*(p;L) := \max\{0,\lambda_V(p;L)-1\}$. 
If $V$ is a cone with apex $p_V$, we put $\lambda_V(p_V;L) = \lambda^*_V(p_V;L) = 0$. 
We clearly have $\lambda_V(p;L) \le \Lambda_V(p)$ and $\lambda^*_V(p;L) \le \Lambda^*_V(p)$, 
for each point $p$.

\begin{lemma}
\label{firstflip} Let $V$ be an irreducible singly ruled
two-dimensional surface of degree $\Deg>1$ in $\reals^3$ or in
$\cplx^3$. Then, for any line $\ell$, except for the (at most) two
exceptional lines of $V$, we have
\begin{align*}
& \sum_{p \in \ell \cap V} \Lambda_V(p) \le \Deg \quad\quad\text{if $\ell$ is not fully contained in $V$} , \\
& \sum_{p \in \ell \cap V} \Lambda^*_V(p) \le \Deg
\quad\quad\text{if $\ell$ is fully contained in $V$} .
\end{align*}
\end{lemma}

\noindent{\bf Proof.} 
To streamline the analysis and avoid degenerate situations that might arise 
over the reals, we confine ourselves to the complex case; as already mentioned,
the incidence bounds that we will obtain will automatically hold over the reals too.
We note that the difference between the two
cases arises because we do not want to count $\ell$ itself---the
former sum would be infinite when $\ell$ is fully contained in $V$.
Note also that if $V$ is a cone and $p_V\in\ell$, we ignore in the
sum the infinitely many lines incident to $p_V$ and contained in $V$.

The proof is a variant of an observation due to Salmon~\cite{salmon}
and repeated in Guth and Katz~\cite{GK2} over the real numbers, and
later in Koll\'ar~\cite{Kollar} over the complex field and other general fields.

By Theorem~\ref{singly}(a), excluding the exceptional lines of $V$,
the set of lines fully contained in $V$ can be parameterized as a
(real or complex) 1-parameter family of generator lines $\ell(t)$,
represented by the irreducible curve $\Sigma_0 \subseteq F(V)$.
Let $V^{(2)}$ denote the locus of points of $V$ that are incident to
at least two generator lines fully contained in $V$. By
Theorem~\ref{singly}(b), $V^{(2)}$ is contained in some
one-dimensional curve $\C \subset V$.

Let $p\in V\cap\ell$ be a point incident to $k$ generator lines of
$V$, other than $\ell$, for some $k\ge 1$. In case $V$ is a cone, we 
assume that $p \ne p_V$.  Denote the generator lines incident to $p$ 
(other than $\ell$, if $\ell\subset V$, in which case it is assumed 
to be a generator) as $\ell_i = \ell(t_i)$, for $t_i\in\Sigma_0$ and 
for $i=1,\ldots,k$. (If $\ell_i$ is a singular point of $F(V)$, it may 
arise as $\ell(t_i)$ for several values of $t_i$, and we pick one arbitrary such value.)
Let $\pi$ be a generic plane containing $\ell$, and consider the
curve $\gamma_0 = V\cap \pi$, which is a plane curve of degree $D$.
Since $V^{(2)}\subseteq\C$ is one-dimensional, a generic choice of
$\pi$ will ensure that $V^{(2)}\cap\pi$ is a discrete set (since
$\ell$ is non-exceptional, it too meets $V^{(2)}$ in a discrete set).

There are two cases to consider: If $\ell$ is fully contained in $V$
(and is thus a generator), then $\gamma_0$ contains $\ell$. In this
case, let $\gamma$ denote the closure of $\gamma_0\setminus\ell$; it
is also a plane algebraic curve, of degree at most $\Deg-1$.
Otherwise, we put $\gamma := \gamma_0$. By Theorem~\ref{singly}(a),
we can take, for each $i=1,\ldots,k$, a sufficiently small open
(real or complex) neighborhood $\Delta_i$ along $\Sigma_0$
containing $t_i$, so that, for any $1 \le i < j \le k$, all the
lines $\ell(t)$, for $t \in \Delta_i \cup \Delta_j$, are distinct.
Put $V_i:= \bigcup_{t\in\Delta_i} \ell(t)$.
Recall that $V_i \cap \pi$ is either a simple arc or a union of
simple arcs meeting at $p$ (depending on whether or not $\ell_i$ is
a regular point of $\Sigma_0$); in the latter case, take $\gamma_i$
to be any one of these arcs. Each of the arcs $\gamma_i$ is incident to
$p$ and is contained in $\gamma$. Moreover, since $\pi$ is
generic, the arcs $\gamma_i$ are all distinct. Indeed, for any $i\ne
j$, and any point $q\in\gamma_i\cap\gamma_j$, there exist $t_i \in
\Delta_i, t_j \in \Delta_j$ such that $\ell(t_i)\cap \pi=
\ell(t_j)\cap \pi = q$, and $\ell(t_i)\ne \ell(t_j)$ (by the
properties of these neighborhoods). Therefore, any point in
$\gamma_i \cap \gamma_j$ is incident to (at least) two distinct
generator lines fully contained in $V$. Again, the generic choice of
$\pi$ ensures that $\gamma_i\cap\gamma_j\subseteq V^{(2)}$ is a
discrete set, so, in particular, $\gamma_i$ and $\gamma_j$ are distinct.

We have therefore shown that (i) if $\ell$ is not contained in $V$
then $p$ is a singular point of $\gamma$ of multiplicity at least
$k$ (for $k\ge 2$; when $k=1$ the point does not have to be
singular), and (ii) if $\ell$ is contained in $V$ then $p$ is
singular of multiplicity at least $k+1$. We have $k \ge
\Lambda_V(p)$ (resp., $k \ge \Lambda^*_V(p)$) if $\ell$ is not fully
contained (resp., is fully contained) in $V$. 
As argued at the beginning of Section~\ref{sec:ruled}, the line $\ell$ can 
intersect $\gamma$ in at most $\Deg$ points, \emph{counted with multiplicity}, 
and the result follows.  \proofend

We also need the following result, established by Guth and
Katz~\cite{GK}; see also \cite{EKS}. It is an immediate consequence
of the Cayley--Salmon--Monge theorem (Theorem~\ref{th:flec2a}) and a
suitable extension of B\'ezout's theorem for intersecting surfaces
(see Fulton~\cite[Proposition 2.3]{Fu84}).
\begin{proposition} [Guth and Katz \protect{\cite{GK2}}]
\label{caysala} Let $V$ be an irreducible two-dimensional variety in
$\cplx^3$ of degree $\Deg$. If $V$ fully contains more than
$11\Deg^2-24\Deg$ lines then $V$ is ruled by (complex) lines.
\end{proposition}

\begin {corollary}
\label{co:caysala} Let $V$ be an irreducible two-dimensional variety
in $\cplx^3$ of degree $\Deg$ that does not contain
any planes. Then the number of lines that are fully contained in the
union of the non-ruled components of $V$ is $O(\Deg^2)$.
\end{corollary}

\noindent{\bf Proof.} Let $V_1,\ldots,V_k$ denote those irreducible
components of $V$ that are not ruled by lines.  By
Proposition~\ref{caysala}, for each $i$, the number of lines fully
contained in $V_i$ is at most $11\deg(V_i)^2-24\deg(V_i)$. Summing
over $i=1,\ldots,k$, the number of lines fully contained in the
union of the non-ruled components of $V$ is at most 
$\sum_{i=1}^k 11\deg(V_i)^2= O(\Deg^2)$. \proofend

The following theorem, which we believe to be of independent
interest in itself, is the main technical ingredient of our analysis.
Note that it holds over both real and complex fields.
\begin{theorem} 
\label{salta} Let $V$ be a possibly reducible two-dimensional
algebraic surface of degree $D>1$ in $\reals^3$ or in $\cplx^3$,
with no linear components. Let $P$ be a set of $m$ distinct points
on $V$ and let $L$ be a set of $n$ distinct lines fully contained in
$V$. Then there exists a subset $L_0\subseteq L$ of at most
$O(\Deg^2)$ lines, such that the number of incidences between $P$
and $L\setminus L_0$ satisfies
\begin{equation} \label{lmlstara}
I(P,L\setminus L_0) = O\left(m^{1/2}n^{1/2}D^{1/2} + m + n\right) .
\end{equation}
\end{theorem}
\noindent{\bf Remark.}
An important feature of the theorem, already noted for the more general Theorem~\ref{th:main3d},
and discussed in more detail later on, is that
the bound in (\ref{lmlstara}) avoids the term $nD$, which arises naturally in
many earlier works, e.g., when bounding the number of incidences between points on $V$
and lines not crossing $V$. This is significant when $D$ is large---see below.

\smallskip

\noindent{\bf Proof.} 
As in the proof of Lemma~\ref{firstflip}, we only work over $\cplx$, 
and the results are then easily transported to the real case too.
Consider the irreducible components
$W_1,\ldots,W_k$ of $V$. By Corollary~\ref{co:caysala}, the number
of lines contained in the union of the non-ruled components of $V$
is $O(\Deg^2)$, and we place all these lines in the exceptional set
$L_0$. In what follows we thus consider only ruled components of
$V$. For simplicity, continue to denote them as $W_1,\ldots,W_k$,
and note that $k\le D/2$.

We further augment $L_0$ as follows. We first dispose of lines of
$L$ that are fully contained in more than one ruled component $W_i$.
We claim that their number is $O(D^2)$. Indeed, for any pair $W_i$,
$W_j$ of distinct components, the intersection $W_i \cap W_j$ is a
curve of degree (at most) $\deg(W_i)\deg(W_j)$, which can therefore
contain at most $\deg(W_i)\deg(W_j)$ lines (by the generalized
version of B\'ezout's theorem~\cite[Proposition 2.3]{Fu84}, already
mentioned in connection with Proposition~\ref{caysala}). Since
$\sum_{i=1}^k \deg(W_i)\le D$, we have
$$
\sum_{i\ne j} \deg(W_i)\deg(W_j) \le \left(\sum_i \deg(W_i)\right)^2
= O(D^2) ,
$$
as claimed. We add to $L_0$ all the $O(D^2)$ lines in $L$ that are
contained in more than one ruled component, and all the exceptional
lines of all singly ruled components. The number of lines of the
latter kind is at most $2k \le 2\cdot(D/2) = D$, so the size of
$|L_0|$ is still $O(D^2)$. Hence, each line of $L_1:=L\setminus L_0$
is fully contained in a unique (singly or doubly) ruled component of
$V$, and is a generator of that component.

The strategy of the proof is to consider each line $\ell$ of $L_1$,
and to estimate the number of its incidences with the points of $P$
in an indirect manner, via Lemma~\ref{firstflip}, applied to $\ell$
and to each of the ruled components $W_j$ of $V$. We recall that
$\ell$ is fully contained in a unique component $W_i$, and treat
that component in a somewhat different manner than the treatment of
the other components.

In more detail, we proceed as follows. We first ignore, for each
singly ruled \emph{conic} component $W_i$, the incidences between
its apex (exceptional point) $p_{W_i}$ and the lines of $L_1$ that
are contained in $W_i$. We refer to these incidences as
\emph{conical incidences} and to the other incidences as
\emph{non-conical}. When we talk about a line $\ell$ incident to
another line $\ell'$ at a point $p$, we will say that $\ell$ is
\emph{conically incident} to $\ell'$ (at $p$) if $p$ is the apex of
some conic component $W_i$ and $\ell'$ is fully contained in $W_i$
(and thus incident to $p$). In all other cases, we will say that
$\ell$ is \emph{non-conically incident} to $\ell'$ (at $p$). (Note
that this definition is asymmetric in $\ell$ and $\ell'$; in
particular, $\ell$ does not have to lie in the cone $W_i$.) We note
that the number of conical point-line incidences is at most $n$,
because each line of $L_1$ is fully contained in a unique component
$W_i$, so it can be involved in at most one conical incidence (at
the apex of $W_i$, when $W_i$ is a cone).

We next prune away points $p\in P$ that are non-conically incident
to at most three lines of $L_1$. (Note that $p$ might be an apex of
some conic component(s) of $V$; in this case $p$ is removed if it is
incident to at most three lines of $L_1$ that are not contained in
any of these components.) We lose $O(m)$ (non-conical) incidences in
this process. Let $P_1$ denote the subset of the remaining points.

\begin{lemma} \label{claim:4D}
Each line $\ell\in L_1$ is non-conically incident, at points of
$P_1$, to at most $4\Deg$ other lines of $L_1$.
\end{lemma}

\noindent{\bf Proof.} Fix a line $\ell\in L_1$ and let $W_i$ denote
the unique ruled component that fully contains $\ell$. Let $W_j$ be
any of the other ruled components. We estimate the number of lines
of $L_1$ that are non-conically incident to $\ell$ and are fully
contained in $W_j$.

If $W_j$ is a regulus, there are at most four such lines, since
$\ell$ meets the quadratic surface $W_j$ in at most two points, each
incident to exactly two generators (and to no other lines contained
in $W_j$). In this case, we write the bound $4$ as $\deg(W_j)+2$.
Assume then that $W_j$ is singly ruled. By Lemma~\ref{firstflip}, we have
$$
\sum_{p \in \ell \cap W_j} \lambda_{W_j}(p;L_1) \le \sum_{p \in \ell
\cap W_j} \Lambda_{W_j}(p) \le \deg(W_j) .
$$
Note that, by definition, the above sum counts only non-conical
incidences (and only with generators of $W_j$, but the exceptional
lines of $W_j$ have been removed from $L_1$ anyway).

We sum this bound over all components $W_j\ne W_i$, including the
reguli. Denoting the number of reguli by $\rho$, which is at most
$D/2$, we obtain a total of
$$
\sum_{j\ne i} \deg(W_j) + 2\rho \le \Deg + 2\rho \le 2\Deg .
$$
Consider next the component $W_i$ containing $\ell$. Assume first
that $W_i$ is a regulus. Each point $p\in P_1\cap\ell$ can be
incident to at most one other line of $L_1$ contained in $W_i$ (the
other generator of $W_i$ through $p$). Since $p$ is in $P_1$, it is
non-conically incident to at least $3-2=1$ other line of $L_1$,
contained in some other ruled component of $V$.
That is, the number of lines that are (non-conically) incident to
$\ell$ and are contained in $W_i$, which apriorily can be
arbitrarily large, is nevertheless at most the number of other lines
(not contained in $W_i$) that are non-conically incident to $\ell$,
which, as shown above, is at most $2\Deg$.

If $W_i$ is not a regulus, Lemma~\ref{firstflip} implies that
$$
\sum_{p \in \ell \cap W_i} \Lambda^*_{W_i}(p) \le \deg(W_i) \le \Deg
,
$$
where again only non-conical incidences are counted in this sum (and
only with generators). That is, the number of lines of $L_1$ that
are non-conically incident to $\ell$ (at points of $P_1$) and are
contained in $W_i$ is at most $\Deg$. Adding the bound for $W_i$,
which has just been shown to be either $D$ or $2D$, to the bound
$2\Deg$ for the other components, the claim follows. \proofend

To proceed, choose a threshold parameter $\xi\ge 3$, to be
determined shortly. Each point $p\in P_1$ that is non-conically
incident to at most $\xi$ lines of $L$ contributes at most $\xi$
(non-conical) incidences, for a total of at most $m\xi$ incidences.
(Recall that the overall number of conical incidences is at most
$n$.) For the remaining non-conical incidences, let $\ell$ be a line in $L_1$
that is incident to $t$ points of $P_1$, so that each such point $p$
is non-conically incident to at least $\xi+1$ lines of $L_1$ (one of
which is $\ell$). It then follows from Lemma~\ref{claim:4D} that
$t\le 4\Deg/\xi$. Hence, summing this over all $\ell \in L_1$, we
obtain a total of at most $4n\Deg/\xi$ incidences. We can now bring
back the removed points of $P\setminus P_1$, since the non-conical
incidences that they are involved in are counted in the bound
$m\xi$. That is, we have
$$
I(P,L_1) \le m\xi + n + \frac{4n\Deg}{\xi} .
$$
We now choose $\xi = (n\Deg/m)^{1/2}$. For this choice to make
sense, we want to have $\xi\ge 3$, which will be the case if $9m\le
n\Deg$. In this case we get the bound
$O\left(m^{1/2}n^{1/2}\Deg^{1/2} + n\right)$. If $9m>n\Deg$ we take
$\xi=3$ and obtain the bound $O(m)$. Combining these bounds, and
adding the at most $n$ conical incidences, the theorem follows.
\proofend

\paragraph{The final stretch: The real case.} 
It remains to bound the number $I(P,L_0)$ of incidences involving the lines in $L_0$. 
We remark that, in both the real and the complex cases, no special properties
need to be assumed for the lines of $L_0$; the only thing that matters is that
their number is small. 
We have $|L_0|=O(\Deg^2)$. 
In the real case, we estimate $I(P,L_0)$ using Guth and Katz's bound (\cite{GK2}; see
Theorem~\ref{ttt}), recalling that no plane contains more than $s$
lines of $L_0$. We thus obtain
\begin{align} \label{lstar}
I(P,L_0) & =
O\left(m^{1/2}|L_0|^{3/4}+m^{2/3}|L_0|^{1/3}s^{1/3}+m+|L_0|\right) \\
& = O\left(m^{1/2}n^{1/2}\Deg^{1/2} +
m^{2/3}\Deg^{2/3}s^{1/3} + m + n \right) \nonumber .
\end{align}
Combining the bounds in Theorem~\ref{salta} and in (\ref{lstar})
yields the asserted bound on $I(P,L)$.

We remark that in the first term in (\ref{lstar}) we have estimated $L_0^{3/4}$
by $O(n^{1/2}D^{1/2})$ instead of the sharper estimate $O(D^{3/2})$. This is
because the term $O(m^{1/2}n^{1/2}D^{1/2})$ appears in Theorem~\ref{salta}
anyway, so the sharper estimate has no effect on the overall asymptotic bound.

\paragraph{The final stretch: The complex case.} 
We next estimate $I(P,L_0)$ in the complex case. Again, we have $n_0:=|L_0|=O(\Deg^2)$. 
We may ignore the points of $P$ that are incident to fewer than three lines of $L_0$,
as they contribute altogether only $O(m)$ incidences. Continue to denote the set of
surviving points as $P$. A point $p$ that is incident to at least three lines of $L_0$ 
is either a singular point of $V$ (when not all its incident lines are coplanar) or 
a flat point of $V$ (see the Appendix, and also Guth and Katz~\cite{GK2}).

We need the following lemma, adapted (with almost the same proof, which we omit) 
from a similar lemma that was established in our 
earlier work \cite{SS4d} for the four-dimensional case (see also Section~\ref{sec:pfmain}).
\begin {lemma}[\protect{\cite[Lemma 2.15]{SS4d}}]
\label {le:lar} Let $f\in \cplx[x,y,z]$ be an irreducible
polynomial. If a line $\ell \subset Z(f)$ is flat, then the tangent
plane $T_p Z(f)$ is fixed for all the non-singular points $p\in\ell$.
\end {lemma}

We decompose $V$ into its irreducible components, and assign each point $p\in P$ 
(resp., line $\ell \in L_0$) to the first component that (fully) 
contains it. Similar to what has been observed above, the number 
of ``cross-incidences'', between points and lines assigned to different 
components of $V$, is $O(n_0 D) = O(D^3)$. We therefore assume, 
as we may, that $V$ is irreducible (over $\cplx$), and write $V=Z(f)$,
for an irreducible trivariate polynomial $f$ of degree $D$.

We call a line \emph{flat} if all its non-singular points are flat.
As argued in the Appendix and in earlier works (see, e.g.,~\cite{EKS,GK2}), since $Z(f)$ is not 
a plane, there exists a certain polynomial $\Pi$ satisfying (i) 
$\deg(\Pi)= 3D-4$, (ii) $Z(f,\Pi)$ is a curve, and (iii) the flat 
points of $P$ and the flat lines of $L_0$ are contained in $Z(f,\Pi)$.

All this implies, arguing as in previous works~\cite{EKS,GK2}, that a
line $\ell \in L_0$ that is non-singular and non-flat contains at
most $4D-4$ points of $P$ (each of which is either singular or flat). We prune away
these lines from $L_0$, losing at most $(4D-4)|L_0| = O(D^3)$ incidences
with the points of $P$. Continue to denote the subset of surviving lines as $L_0$.

Therefore, it remains to bound the number of incidences between the surviving
points and the surviving lines, each of which is either singular or flat. Write $P$ 
as the union of the subset $P_f$ of flat points and the subset $P_s$ of singular points. 
Similarly, write $L_0$ as the union of the subset $L_f$ of flat lines and the subset $L_s$ 
of singular lines. A singular line contains no flat points, and a flat
line contains at most $D-1$ singular points. Thus, 
$$
I(P,L_0) \le I(P_f, L_f) + I(P_s, L_s) + n_0 D .
$$
By Lemma~\ref{le:lar}, all the non-singular points of a flat line
have the same tangent plane. Assign each point $p\in P_f$ (resp.,
line in $L_f$) to its tangent plane $T_p Z(f)$ (resp., $T_p Z(f)$
for some non-singular point $p\in P_f \cap \ell$; we only consider
lines in $L_f$ that are incident to at least one point in $P_f$).
We have therefore partitioned the points in $P_f$ and the lines in
$L_f$ among planes in some finite set $H=\{h_1, \ldots, h_k\}$, and we
only need to count incidences between points and lines assigned to
the same plane. Within each $h\in H$, we have a set $P_h\subseteq P_f$ 
of $m_h$ points in $h$, and a set $L_h\subseteq L_f$ of $n_h$ lines 
contained in $h$. Using the planar bound (\ref{inc2}), which also holds in the complex plane,
the number of incidences within $h$ is $O(m_h^{2/3}n_h^{2/3}+m_h+n_h)$.  
Summing these bounds over $h\in H$, and using H\"older's inequality, and the fact that
$n_h \le s$ for each $h$, we obtain a total of $O(m^{2/3}n_0^{1/3}s^{1/3}+m+n_0)$ incidences.

\paragraph{Bounding incidences involving singular points and lines.}
Bounding the number of incidences between the singular points and
lines is done via degree reduction. Assuming, without loss of generality, that
$f_x$ does not vanish identically on $Z(f)$, the points of $P_s$ and the
lines of $L_s$ are then (fully) contained in $Z(f_x)$, and $\deg(f_x) \le D-1$. 
We thus construct a sequence of partial derivatives of $f$ that are not 
identically zero on $Z(f)$. For this we assume, as we may, that $f$, and 
each of its derivatives, are square-free; whenever this fails, we replace the 
corresponding derivative by its square-free counterpart before continuing to 
differentiate. Without loss of generality, assume that this sequence is
$f,f_x,f_{xx}$, and so on.  Denote the $j$-th element in this
sequence as $f_j$, for $j=0,1,\ldots$ (so $f_0=f$, $f_1=f_x$, and so
on). Assign each point $p\in P$ to the first polynomial $f_j$ in the
sequence for which $p$ is non-singular; more precisely, we assign
$p$ to the first $f_j$ for which $f_j(p)=0$ but $f_{j+1}(p)\ne 0$
(recall that $f_0(p)$ is always $0$ by assumption. Similarly, assign
each line $\ell$ to the first polynomial $f_j$ in the sequence for
which $\ell$ is fully contained in $Z(f_j)$ but not fully contained
in $Z(f_{j+1})$ (again, by assumption, there always exists such a polynomial
$f_j$). If $\ell$ is assigned to $f_j$ then it can only contain
points $p$ that were assigned to some $f_k$ with $k\ge j$. Indeed,
if $\ell$ contained a point $p$ assigned to $f_k$ with $k<j$ then
$f_{k+1}(p)\ne 0$ but $\ell$ is fully contained in $Z(f_{k+1})$,
since $k+1\le j$; this is a contradiction that establishes the claim.

Fix a line $\ell\in L$, which is assigned to some $f_j$. An incidence 
between $\ell$ and a point $p\in P$, assigned to some $f_k$, for $k>j$, 
can be charged to the intersection of $\ell$ with $Z(f_{j+1})$ at $p$ 
(by construction, $p$ belongs to $Z(f_{j+1})$). The number of such 
intersections is at most $\deg(f_{j+1}) \le D-j-1 \le D$, so the overall 
number of incidences of this sort, over all lines $\ell\in L$, is
$O(n_0D) = O(D^3)$. It therefore suffices to consider only incidences between
points and lines that are assigned to the same zero set $Z(f_i)$.

The reductions so far have produced a finite collection of up to
$D$ polynomials, each of degree at most $D$, so that the points
of $P$ are \emph{partitioned} among the polynomials and so are the
lines of $L$, and we only need to bound the number of incidences
between points and lines assigned to the \emph{same} polynomial.
Moreover, for each $j$, all the points assigned to $f_j$ are non-singular, by construction.
For each $j$, let $P_j$ and $L_j$ denote the subsets of $P$ and of $L_0$, respectively,
that are assigned to $f_j$, and put $m_j := |P_j|$ and $n_j := |L_j|$.
We have $\sum_j m_j \le m$ and $\sum_j n_j \le n_0$.

We would like to apply the preceding analysis to $P_j$ and $L_j$, but we face
the technical issue that $f_j$ might be reducible and have some linear factors.
(The theorem does not require the variety to be irreducible, but forbids it to have linear components.)

We therefore proceed as follows. We first consider only those points and lines that are (fully)
contained in some nonlinear component of $Z(f_j)$. We apply the preceding analysis to 
these sets, and obtain the incidence bound 
$$
O\left( m_j^{2/3}n_j^{1/3}s^{1/3} + m_j + n_jD \right) .
$$
Summing these bounds over all $f_j$'s, using H\"older's inequality, we get the overall bound
$$
O\left( m^{2/3}n_0^{1/3}s^{1/3} + m + n_0D \right) .
$$
To bound the number of incidences involving points and lines in the linear components of $Z(f_j)$,
for any fixed $j$, we reason as we did above, when handling the flat points and lines.
That is, we order arbitrarily the planar components of $Z(f_j)$, assign each point to 
first component that contains it, and assign each line to the first component that fully contains
it. As $Z(f_j)$ has at most $D$ such components, the number of ``cross incidences,'' between 
points and lines assigned to different components, is $O(n_jD)$. 
For the number of ``same component'' incidences, we use the T\'oth-Zahl extension of the
Szemer\'edi-Trotter bound~\cite{Tot,Zah} in each plane, and sum them up, exactly as in 
the case of flat points and lines, and get a total bound of
$$
O\left( m_j^{2/3}n_j^{1/3}s^{1/3} + m_j + n_jD \right) ,
$$
and, summing these bounds over all $f_j$'s, as above, we get the overall bound
$$
I(P_s,L_s) = O\left( m^{2/3}n_0^{1/3}s^{1/3} + m + n_0D \right) ,
$$
thereby completing the proof of part (b) (for the case where the ambient space is three-dimensional).


\paragraph{Reduction to three dimensions.} 
To complete the analysis, we need to consider the case where $V$ is a
two-dimensional variety embedded in $\reals^d$, for $d>3$.

Let $H$ be a generic 3-flat, and denote by $P^*, L^*$, and $V^*$ the
respective projections of $P, L$, and $V$ onto $H$. Since $H$ is
generic, we may assume that all the projected points in $P^*$ are
distinct, and so are all the projected lines in $L^*$. Clearly,
every incidence between a point of $P$ and a line of $L$ corresponds
to an incidence between the projected point and line. Since no
2-flat contains more than $s$ lines of $L$, and $H$ is generic,
repeated applications of Lemma~\ref{ap:th} imply that no plane in
$H$ contains more than $s$ lines of $L^*$.

One subtle point is that the set-theoretic projection $V^*$ of $V$
does not have to be a real algebraic variety (in general, it is only
a semi-algebraic set), but it is always contained in a
two-dimensional real algebraic variety $\tilde V$, which we call, as
we did in an earlier work~\cite{SSsocg}, the \emph{algebraic
projection} of $V$; it is the zero set of all polynomials belonging
to the ideal of polynomials vanishing on $V$, after eliminating
variables in the complementary space of $H$ (this is also known as
an \emph{elimination ideal} of $V$; see Cox et al.~\cite{CLO} for
details), and is equal to the Zariski closure of $V^*$. 
Since the closure of a projection does not increase the
original degree (see, e.g., Harris~\cite{Har}), $\deg(\tilde V)\le D$. 
That $\tilde V$ does not contain a 2-flat follows by a suitable
adaptation of the argument in~\cite[Lemma 2.1]{SSsocg}
(which is stated there for $d=4$ over the reals), that applies for
general $d$ and over the complex field too.

In conclusion, we have $I(P,L) \le I(P^*,L^*)$, where $P^*$ is a set
of $m$ points and $L^*$ is a set of $n$ lines, all contained in the
two-dimensional algebraic variety $\tilde V$, embedded in 3-space,
which is of degree at most $D$ and does not contain any plane, and
no plane contains more than $s$ lines of $L^*$. The preceding
analysis thus implies that the bound asserted in the theorem applies
in any dimension $d\ge 3$. \proofend

\section{Proof of Theorem~\ref{th:main4d}} \label{sec:pfmain}

In most of this section we assume that the ambient space is
four-dimensional, and work, with a few exceptions, over the complex field.
The reduction from higher dimensions to four dimensions is handled as the 
reduction to three dimensions just discussed. 

We exploit the following useful corollary of
Theorem~\ref{th:main3d} (recall that we are now in four dimensions).
Note that in the bound given below, the term $nD$ does not appear yet.

\begin {corollary}
\label{th:degsquared} 
Let $f$ and $g$ be two $4$-variate polynomials, over $\reals$ or $\cplx$,
of degree $O(D)$, such that 
$Z(f,g)$ is two-dimensional over $\cplx$. Let $P$ be a set of $m$ points and $L$ a set of $n$ lines,
such that all the points of $P$ and all the lines of $L$ are (fully) contained in 
the union of the irreducible components of $Z(f,g)$ that are not 2-flats.
Assume also that no 2-flat contains more than $s$ lines of $L$. Then we have \\
{\bf (a) in the real case:}
\begin {equation} \label{cor:3d}
I(P,L) = O\left(m^{1/2}n^{1/2}D+m^{2/3}D^{4/3}s^{1/3}+m+n\right) ,
\end {equation}
{\bf (b) and in the complex case:}
\begin {equation} \label{cor:3dcplx}
I(P,L) = O\left(m^{1/2}n^{1/2}D+m^{2/3}D^{4/3}s^{1/3}+D^6+m+n\right).
\end {equation}
\end {corollary}

\noindent{\bf Proof.}
Let $Z(f,g) = \bigcup_{i=1}^s V_i$ be the
decomposition of $Z(f,g)$ into its irreducible components. By the
generalized version of B\'ezout's theorem~\cite{Fu84}, we have
$\sum_{i=1}^s \deg(V_i) \le \deg(f)\deg(g)=O(D^2)$. Assume that
$V_1,\ldots, V_k$ are the components that are not 2-flats, for some $k\le s$,
and let $W$ denote their union. As just observed, $\deg(W) = O(D^2)$.
Applying to $W$ Theorem~\ref{th:main3d}(a) (over the reals) or
Theorem~\ref{th:main3d}(b) (over the complex) thus completes the proof.
\proofend

\smallskip

\noindent{\bf Remark.}
Corollary~\ref{th:degsquared}(a) is significant when $D\ll n^{1/4}$. For larger
values of $D$, we can project $P$, $L$, and $Z(f,g)$ onto some generic 3-flat,
and apply the incidence bound of Guth and Katz, as given in Theorem~\ref{ttt},
within that 3-flat. When $D\ge n^{1/4}$, the resulting bound is better than
the one in (\ref{cor:3d}).

The proof of Theorem~\ref{th:main4d} revisits the proof of
Theorem~\ref{th:mainfocs}, as presented in~\cite{SS4d}, and applies
Corollary~\ref{th:degsquared} as a major technical tool. The only difference 
between the real and the complex cases is in the application of that corollary. 
Except for this application, we will work over the complex domain, but the
analysis carries over, in a straightforward manner, to the real case too.

Each line $\ell \in L$ that is not fully contained in $V$ contributes at most
$D$ incidences, for a total of $O(nD)$ incidences. We thus assume,
as we may, that all the lines of $L$ are contained in $V$. Let
$V=\bigcup_{i=1}^t V_i$ be the decomposition of $V$ into its
irreducible components, and assign each point $p\in P$, (resp., line
$\ell \in L$) to the first $V_i$ that contains it (resp., fully
contains it; such a $V_i$ always exists). It is easy to verify that
points and lines that are assigned to different $V_i$'s contribute
at most $nD$ incidences. Indeed, any such incidence $(p,\ell)$ can
be charged to an intersection point of $\ell$ with the component
$V_i$ that $p$ is assigned to, and thus there are at most $\sum_i \deg(V_i) = D$ 
such incidences for each $\ell$, for an overall number
of $O(nD)$ such incidences. Therefore, it suffices to establish the
bound in (\ref{ma:in}) or in (\ref{ma:incplx}) for the number of incidences between points
and lines assigned to the same component. The exponents in the first two terms in
either of these bounds are favorable, in the sense that they sum up to $1$
(ignoring the factors $D$ and $s^{1/3}$), thereby allowing us to estimate
the sum of the resulting bounds, over the components of $V$, via
H\"older's inequality, to obtain the bounds in Theorem~\ref{th:main4d}
(see below for details). We thus assume that $V$ is irreducible, and write 
$V=Z(f)$, for some real or complex irreducible polynomial $f$ of degree $D$.

We assume for now that $P$ consists exclusively of
\emph{non-singular} points of the \emph{irreducible} variety $Z(f)$.
The treatment of the singular points, similar to the handling of singular 
points in the proof of Theorem~\ref{th:main3d}(b), will be given towards 
the end of the proof.

We recall the defintion of the \emph{four-dimensional flecnode polynomial} 
$\fl{f}^4$ of $f$, as given in Section~\ref{Se:pre}. That is,
$\fl{f}^4$ vanishes at each \emph{flecnode} $p\in Z(f)$, namely, points 
$p$ for which there exists a line that is incident to $p$ and osculates 
to $Z(f)$ up to the fourth order. It is obtained by eliminating (in a 
standard manner, using resultants) the direction $v$ of the osculating 
line from the four homogeneous equations given by the vanishing of the 
first four terms of the Taylor expansion of $f(p+tv)$ near $p$ ($v$ is 
a point in projective $3$-space, and its elimination from these four 
equations yields a single polynomial equation in $p$). Clearly, $\fl{f}^4$ 
vanishes identically on every line of $L$, and thus also on $P$ (assuming 
that each point of $P$ is incident to at least one line of $L$).
As in the three-dimensional case, its degree can be shown to be $O(D)$.

If $\fl{f}^4$ does not vanish identically on $Z(f)$, then
$Z(f,\fl{f}^4) := Z(f)\cap Z(\fl{f}^4)$ is a two-dimensional variety
that contains $P$ and all the lines of $L$, and is of degree
$\deg(f)\cdot\deg(\fl{f}^4) = O(\Deg^2)$ (so we are in the setup assumed
in Corollary~\ref{th:degsquared}). The other possibility is
that $\fl{f}^4$ vanishes identically on $Z(f)$, and then a theorem
of Landsberg~\cite{Land} (see also~\cite{SS4d} for details) implies that
$Z(f)$ is ruled by (real or complex) lines. (Landsberg's theorem is a
generalization of Theorem~\ref{th:flec2a}, the classical Cayley--Salmon--Monge
theorem~\cite{Mon,salmon} in three dimensions.)

\smallskip

\noindent{\bf First case: $Z(f,\fl{f}^4)$ is two-dimensional.}
Put $g=\fl{f}^4$ and apply Corollary~\ref{th:degsquared} to $f$ and $g$. 
In the real case we obtain the bound
\begin {equation*}
O\left( m^{1/2}n^{1/2}D+m^{2/3}D^{4/3}s^{1/3}+m+n \right) ,
\end {equation*}
and in the complex case we obtain the bound
\begin {equation*}
O\left( m^{1/2}n^{1/2}D+m^{2/3}D^{4/3}s^{1/3}+D^6+m+n \right) ,
\end {equation*}
over all components of $Z(f,\fl{f}^4)$ that are not 2-flats.

\smallskip

\noindent{\bf Incidences within 2-flats fully contained in $Z(f,\fl{f}^4)$.}
From this point on in the proof, there is no distinction between the 
real and complex cases, so we work over the complex domain.
The strategy here is to distribute
the points of $P$ and the lines of $L$ among the 2-flats that
contain them (lines not contained in any 2-flat are fully contained in
some other component of $Z(f,\fl{f}^4)$ and are dealt with, as above, within that component).
See also below for a more detailed account of this strategy.
Points that belong to at most two such 2-flats get duplicated at
most twice, and we bound the number of incidences with these points by applying
the planar bound (\ref{inc2}) (which, as we recall, also holds in the complex plane~\cite{Tot,Zah})
to each 2-flat separately, and sum up the bounds, to get
$O(m^{2/3}n^{1/3}s^{1/3} + m + n)$, using H\"older's inequality, combined with the
assumption that no 2-flat contains more than $s$ lines of $L$.

Extending to four dimensions (and in the complex domain) the notation for the 
three-dimensional case from Guth and Katz~\cite{GK} (see the appendix, 
\cite{EKS}, and also Pressley~\cite{Pr} and Ivey and Landsberg~\cite{IL} 
for more basic references), we call a non-singular point $p$ of $Z(f)$ 
\emph{linearly flat}, if it is incident to at least three distinct 2-flats 
that are fully contained in $Z(f)$ (and thus also in the tangent hyperplane 
$T_p Z(f)$).  Linearly flat points can then be shown to be \emph{flat}, meaning
that the second fundamental form of $f$ vanishes at them (see the appendix
and~\cite{SS4d}). This property, at a point $p$, can be expressed by several 
polynomials of degree $3D-4$ vanishing at $p$ (see~\cite[Section 2.5]{SS4d}).
We call a line \emph{flat} if all its non-singular points are flat.
Each line of $L$ that is not flat contains at most $O(D)$ flat
points, and thus the non-flat lines contribute a total of at most
$O(nD)$ incidences with flat points, so we assume in what follows
that the points of $P$ and the lines of $L$ are all flat. Since
$Z(f)$ is not a hyperplane, the second fundamental form does not
vanish identically on $Z(f)$ (this property holds in any dimension;
see, e.g.~\cite[Exercise 3.2.12.2]{IL}), and it then follows from the
characterization of flat points that there exists a certain
polynomial $\Pi$ satisfying (i) $\deg(\Pi)= 3D-4$, (ii) $Z(f,\Pi)$
is two-dimensional, and (iii) the (flat) points of $P$ and the
(flat) lines of $L$ are contained in $Z(f,\Pi)$.

In this case we partition the points of $P$ among their tangent hyperplanes, 
and denote the resulting set of hyperplanes by $H$. Similar to Lemma~\ref{le:lar}, 
we have the property that all the (flat) points that lie on the same flat line 
have the \emph{same} tangent hyperplane to $Z(f)$ (see~\cite[Lemma 2.15]{SS4d}).
We use this property to obtain a partition of the points and lines into 
distinct (tangent) hyperplanes, so that it suffices
to bound the number of incidences within each hyperplane in $H$.

For each $h \in H$, we have a set $P_h\subseteq P$ of $m_h$ points
in $h$, a set $L_h\subseteq L$ of $n_h$ lines contained in $h$, and
a set $F_h$ of 2-flats contained in $Z(f,\fl{f}^4) \cap h$; the
2-flats of $F_h$ are the 2-flat components of $Z(f,\fl{f}^4) \cap h$
(ignoring singularities, which are treated later, these 2-flats are necessarily contained in
the corresponding tangent hyperplanes). Notice that each 2-flat in
$F_h$ is also contained in the two-dimensional surface $Z(f)\cap h$,
which is of degree $D$, so, by the generalized version of B\'ezout's
theorem~\cite{Fu84}, we have $|F_h|\le D$. We assign each point $p
\in P_h$ (resp., line $\ell \in L_h$) to the first 2-flat in $F_h$
that (fully) contains it. Similar to what has been observed above,
the number of ``cross-incidences'', between points and lines
assigned to different $2$-flats, within $h$ is at most $n_h D$, for
a total, over the hyperplanes $h\in H$, of at most $nD$ incidences.
Again, using the planar bound (\ref{inc2}) and H\"older's
inequality, the numbers of incidences within the 2-flats of $h$ sum
up to $O(m_h^{2/3}n_h^{2/3}+m_h+n_hD)$, and, summing over the
hyperplanes $h\in H$, using H\"older's inequality once again, and
the fact that $n_h \le s$, we obtain a total of
$O(m^{2/3}n^{1/3}s^{1/3}+m+nD)$ incidences.

\smallskip

\noindent{\bf Remark.}
The novel feature of this step of the proof, as compared with the
analogous argument used in \cite{SS4d}, is that the number of 2-flats in
$F_h$, for any fixed $h$, is at most $D$. This allows us to
bound the number of incidences within each hyperplane $h\in H$ separately,
so that, within each such hyperplane, instead of using the Guth-Katz
bound (in the real case), we partition the points and lines among at most 
$D$ planes, and then use the Szemer\'edi-Trotter bound in the real case,
or the T\'oth-Zahl bound in the complex case.
The fact that there are at most $D$ planes within each hyperplane $h\in H$ 
guarantees that the number of ``cross-incidences'' (within $h$) is at most 
$n_h D$, for a total of $nD$ incidences (notice that, in contrast, the 
\emph{total} number of planes (over all $h\in H$) can be arbitrarily large).

\smallskip

\noindent{\bf Second case: $Z(f)$ is ruled by (complex) lines.} We
next consider the case where the four-dimensional flecnode
polynomial $\fl{f}^4$ vanishes identically on $Z(f)$. By Landsberg's
theorem~\cite{Land} mentioned above, this
implies that $Z(f)$ is ruled by (complex) lines.

As in the three-dimensional case treated in the appendix (see, e.g., the proof of Theorem~\ref{singly}), 
we denote by $\Sigma_p^3$ (resp., $\Sigma_p$), for $p\in Z(f)$, the set of all lines that 
are incident to $p$ and osculate to $Z(f)$ to order $3$ at $p$ (resp., are 
contained in $Z(f)$). We put $\Sigma^3:= \bigcup_{p\in Z(f)} \Sigma_p^3$, and
$\Sigma:= \bigcup_{p\in Z(f)} \Sigma_p$. ($\Sigma$ is the \emph{Fano variety} 
of (lines contained in) $Z(f)$, now represented in a higher-dimensional projective space.)

In~\cite{SS4d}, we proved that, for each $p\in P$, either $|\Sigma_p|\le 6$ or 
$\Sigma_p$ is infinite.  In the interest of completeness, we recall here the 
outline of this argument. The analysis provides an algebraic characterization
of points $p$ of the latter kind, which uses an auxiliary polynomial
$U=U(p;u_0,u_1,u_2,u_3)$, called the \emph{$u$-resultant}, defined in terms of $f$
and its derivatives at $p$ (see~\cite{SS4d} and also~\cite{CLO} for details), where
$(u_0,\ldots,u_3)$ denotes the direction of a line incident to $p$ (in homogeneous
coordinates). The polynomial $U$ is of degree $O(D)$ in $p$ and is a homogeneous
polynomial of degree six in $u$. The characterization is that $\Sigma_p^3$ is
infinite if and only if $U(p;u_0,u_1,u_2,u_3)\equiv 0$, as a polynomial of $u$, at
$p$. In the complementary case, B\'ezout's theorem~\cite{Fu84} can be used to
show that there are only six lines in $\Sigma_p^3$, and thus at most
six lines in $\Sigma_p$. Pruning away points $p \in P$ with
$|\Sigma_p^3|\le 6$ (the number of incidences involving these points
is at most $6m = O(m)$), we may then assume that $\Sigma_p^3$ is
infinite for every $p\in P$.

If $U(p;u_0,u_1,u_2,u_3)$ does not vanish identically (as a
polynomial in $u_0,u_1,u_2,u_3$) at every point $p \in Z(f)$, then
at least one of its coefficients, call it $c_U$, which is a
polynomial in $p$, of degree $O(D)$, does not vanish identically on
$Z(f)$. In this case, as $U$ vanishes identically at every point of
$P$ (as a polynomial in $u_0,u_1,u_2,u_3$), we have $P \subset
Z(f,c_U)$, which is a two-dimensional variety. The machinery
developed in the first case can then be applied here (with $g=c_U$),
and the bounds and properties derived for that case hold here too.

We may therefore assume that $U(p;u_0,u_1,u_2,u_3) \equiv 0$ at
every non-singular point $p \in Z(f)$ (as a polynomial in
$u_0,u_1,u_2,u_3$). By the aforementioned characterization via
$u$-resultants, it follows that $\Sigma_p^3$ is then infinite at each such point.

We now use another theorem of Landsberg~\cite[Theorem 3.8.7]{IL}: Let $f$ be 
a polynomial over $\mathbb P^{4}(\cplx)$, such that there exists an irreducible
component $\Sigma_{0}^3 \subset \Sigma^3 = \Sigma^3(Z(f))$ with the
property that, for every point $p$ in a Zariski-open set\footnote{%
  See the appendix and Cox et al.~\cite[Section 4.2]{CLO} for further details.}
$\mathcal O \subset Z(f)$, $\dim \Sigma_{0,p}^3 \ge 1$, where $\Sigma_{0,p}^3$
is the set of lines in $\Sigma_0^3$ incident to $p$. Landsberg's
theorem then asserts that, for every point $p\in \mathcal O$, all
lines in $\Sigma_{0,p}^3$ are fully contained in $Z(f)$; that is, for each $p\in {\mathcal O}$,
$\Sigma_{0,p}^3$ is equal to the set $\Sigma_{0,p}$ of lines incident to $p$ and fully contained in $Z(f)$.

Since $\Sigma_p^3$ is infinite at each non-singular point $p\in Z(f)$, 
its dimension is $\ge 1$ at each such point. As shown in~\cite{SS4d}, the 
main condition in Landsberg's theorem, about the existence of a component 
$\Sigma_0^3$ of $\Sigma^3$ with the required properties, is satisfied too. 
One can then argue that the conclusion of Landsberg's theorem holds at \emph{every} 
point of $Z(f)$; see Lemma~\ref{le:rs} in the appendix for a similar claim 
concerning two-dimensional surfaces.  That is, $Z(f)$ is \emph{infinitely ruled} 
by (complex) lines, in the sense that each point $p \in Z(f)$ is incident to 
infinitely many (complex) lines that are fully contained in $Z(f)$, and, moreover, 
$\Sigma_{0,p}^3 = \Sigma_{0,p}$ at each $p$. That is, $\Sigma_0^3$ is contained in 
$\Sigma$. Denoting this set as $\Sigma_0$, it is shown (in full detail) in~\cite{SS4d} 
that the union of the lines in $\Sigma_0$ is equal to $Z(f)$, and that $\dim(\Sigma_0) \ge 3$.


\smallskip

\noindent{\bf Severi's theorem.} The following theorem was already
used in~\cite{SS4d}, and we make a similar use thereof
here too. It has been obtained by Severi \cite{severi} in 1901. A
variant of this result has also been obtained by Segre~\cite{Seg};
see also the more recent works \cite{MP,Ric,Ro}.
We state here a special case of the theorem that we need.

\begin {theorem}[Severi's Theorem~\protect{\cite{severi}}; special case]
\label{th:sev} Let $X \subset \mathbb P^4(\cplx)$ be a three-dimensional 
irreducible variety, and let $\Sigma_0$ be a component of maximal dimension 
of the Fano variety $\Sigma=\Sigma(X)$ of $X$, such that the lines of 
$\Sigma_0$ cover $X$. Then the following holds. (i) If $\dim (\Sigma_0) = 4$, 
then $X$ is a hyperplane. (ii) If $\dim(\Sigma_0) = 3$, then either $X$ is 
a quadric, or $X$ is ruled by 2-flats.
\end {theorem}

Informally, $\dim(\Sigma_0) = 3$ corresponds to the case where $X$
is infinitely ruled by lines of $\Sigma_0$: There are four degrees
of freedom to specify a line in $\Sigma_0$, three to specify $p\in X$,
and one to specify the line in $\Sigma_{0,p}$. (We can assume
that $\Sigma_{0,p}$ is one-dimensional, because if it were
two-dimensional, then $X$ would have been a hyperplane.) However,
one degree of freedom has to be removed, to account for the fact
that the same line (being contained in $X$) arises at each of its points. 
Severi's theorem
asserts, again informally, that in this case the infinite family of
lines of $\Sigma_{0,p}$ must form a 2-flat, unless $X$ is a quadric
or a hyperplane. Moreover, by~\cite[Theorem 3.9]{SS4d} (whose proof
is based on Theorem~\ref{th:harr}, given in the appendix),
$\Sigma_0$ has maximal dimension.

Applying the second case in Severi's theorem to $Z(f)$, which is
justified by the preceding arguments, we conclude that either $Z(f)$
is a quadric or it is ruled by 2-flats. The cases where $Z(f)$ is a
quadric or a hyperplane are ruled out by our assumption, so we
only need to consider the case where $Z(f)$ is ruled by (complex) $2$-flats.

\paragraph{The case where $Z(f)$ is ruled by $2$-flats.} \label{2flat}
Handling this last step is somewhat intricate; it resembles the
analysis of flat points and lines in the first case, where here
points and lines are partitioned among the ruling 2-flats. In this
case, every point $p \in Z(f)$ is incident to at least one 2-flat
$\tau_p\subset Z(f)$. Let $D_p$ denote the set of 2-flats that are
incident to $p$ and are contained in $Z(f)$.

For a non-singular point $p\in Z(f)$, if $\vert D_p \vert > 2$, then
$p$ is a (linearly flat and thus) flat point of $Z(f)$. Recall that
we have bounded the number of incidences involving flat points (and
lines) by partitioning them among a finite number of containing hyperplanes,
and by bounding the incidences within each hyperplane. Lines
incident to fewer than $3D-4$ points of $P$ have been pruned away,
losing only $O(nD)$ incidences, and the remaining lines are all flat. 
Repeating this argument here, noticing that here too, the number of 
2-flats contained in a hyperplane is at most $D$, we obtain the bound
$$
O\left(m^{2/3}n^{1/3}s^{1/3} + m + nD\right).
$$
In what follows we therefore assume that all points of $P$ are non-singular 
and non-flat (call these points, as in \cite{SS4d}, \emph{ordinary} for
short), and therefore $\vert D_p \vert = 1$ or $2$, for each such
$p$. Put $H_1(p)$ (resp., $H_1(p), H_2(p)$) for the 2-flat (resp.,
two 2-flats) in $D_p$, when $\vert D_p \vert = 1$ (resp., $\vert D_p \vert = 2$).

Clearly, each line in $L$, containing at least one ordinary point $p
\in Z(f)$, is fully contained in at most two 2-flats fully contained
in $Z(f)$ (namely, the 2-flats of $D_p$).

Assign each ordinary point $p\in P$ to \emph{each} of the at most
two 2-flats in $D_p$, and assign each line $\ell \in L$ that is
incident to at least one ordinary point to the at most two 2-flats
that fully contain $\ell$ and are fully contained in $Z(f)$ (it is
possible that $\ell$ is not assigned to any $2$-flat---see below).
Changing the notation, enumerate these 2-flats, over all ordinary
points $p\in P$, as $U_1,\ldots,U_k$, and, for each $i=1,\ldots,k$,
let $P_i$ and $L_i$ denote the respective subsets of points and
lines assigned to $U_i$, and let $m_i$ and $n_i$ denote their
cardinalities. We then have $\sum_i m_i \le 2m$ and $\sum_i n_i \le
2n$, and the total number of incidences within the 2-flats $U_i$
(excluding lines not assigned to any 2-flat) is at most
$\sum_{i=1}^k I(P_i,L_i)$. This incidence count can be obtained
exactly as in the first case of the analysis, with the aid of
H\"older's inequality, and yields the bound
$$
\sum_{i=1}^k I(P_i,L_i) = O\left( m^{2/3}n^{1/3}s^{1/3} + m + n \right) .
$$
As noted, this bound does not take into account incidences involving
lines which are not contained in any of the 2-flats $U_i$ (and are
therefore not assigned to any such 2-flat). It suffices to consider
only lines of this sort that are non-singular and non-flat, since
singular or flat lines are only incident to singular or flat points,
and we assumed above that all the points of $P$ are ordinary points.
If $\ell$ is a non-singular and non-flat line, and is not fully
contained in any of the $U_i$, we call it a \emph{piercing line} of
$Z(f)$. We need the following lemma from~\cite{SS4d}.

\begin {lemma} [\protect{\cite[Lemma 3.13]{SS4d}}]
\label{le:pierc2} Let $p\in Z(f)$ be an ordinary point. Then $p$ is
incident to at most one piercing line.
\end {lemma}

Therefore, each ordinary point $p\in P$ is incident to at most one
piercing line, and the total contribution of incidences involving
ordinary points and piercing lines is at most $m$.


In conclusion, combining the bounds that we have obtained for the
various subcases of the second case, we get the desired bound in (\ref{ma:in}).

\paragraph{Incidences involving singular points of $Z(f)$.} \label{sing}
The forthcoming reasoning is very similar to the handling of singular points and lines
in the proof of Theorem~\ref{th:main3d}(b), although it is somewhat more involved
because we need to ensure that the resulting polynomials that we construct be irreducible;
we present the analysis in detail, for the sake of clarity.

In the analysis presented so far, we have assumed that the points
of $P$ are non-singular points of $Z(f)$. To reduce the general
setup to this situation we proceed as follows; an identical
reduction has also been used in~\cite{SS4d}. We only handle lines that are fully
contained in $Z(f)$, because the other lines contribute at most $O(nD)$ incidences.
We construct a sequence of partial derivatives of $f$ that are not identically zero on
$Z(f)$. For this we assume, as we may, that $f$, and each of its
derivatives, are square-free; whenever this fails, we replace the
corresponding derivative by its square-free counterpart before
continuing to differentiate.
Without loss of generality, assume that this sequence is obtained by
always differentiating with respect to $x$, and denote
the $j$-th element in this sequence as $f_j$, for $j=0,1,\ldots$.
That is, $f_0=f$, $f_1=f_x$, with repeated factors removed,
$f_2$ is the $x$-derivative of $f_1$, again with repeated factors
removed, and so on. Assign each point $p\in P$ to
the first polynomial $f_j$ in the sequence for which $p$ is
non-singular; more precisely, we assign $p$ to the first $f_j$ for
which $f_j(p)=0$ but $f_{j+1}(p)\ne 0$ (recall that $f_0(p)$ is
always $0$ by assumption). Similarly, assign each line $\ell$ to the
first polynomial $f_j$ in the sequence for which $\ell$ is fully
contained in $Z(f_j)$ but not fully contained in $Z(f_{j+1})$
(again, by assumption, there always exists such a $j$). If $\ell$ is
assigned to $f_j$ then it can only contain points $p$ that were
assigned to some $f_k$ with $k\ge j$. Indeed, if $\ell$ contained a
point $p$ assigned to $f_k$ with $k<j$ then $f_{k+1}(p)\ne 0$ but
$\ell$ is fully contained in $Z(f_{k+1})$, since $k+1\le j$; this is
a contradiction that establishes the claim.

Fix a line $\ell\in L$, which is assigned to some $f_j$. An
incidence between $\ell$ and a point $p\in P$, assigned to some
$f_k$, for $k>j$, can be charged to the intersection of $\ell$ with
$Z(f_{j+1})$ at $p$ (by construction, $p$ belongs to $Z(f_{j+1})$).
The number of such intersections is at most $\deg(f_{j+1}) \le D-j-1 \le D$, 
so the overall number of incidences of this sort, over all lines $\ell\in L$, is
$O(nD)$. It therefore suffices to consider only incidences between
points and lines that are assigned to the same zero set $Z(f_i)$.

The reductions so far have produced a finite collection of up to
$O(D)$ polynomials, each of degree at most $D$, so that the points
of $P$ are \emph{partitioned} among the polynomials and so are the
lines of $L$, and each point $p$ is non-singular with respect to the
polynomial it is assigned to, and we only need to bound the number of incidences
between points and lines assigned to the \emph{same} polynomial.
This is not the end yet, because the various (reduced forms of the)
partial derivatives might be reducible, which we want to avoid. Thus, in a final
decomposition step, we split each derivative polynomial $f_j$ into
its irreducible factors, and reassign the points and lines that were
assigned to $Z(f_j)$ to the various factors, by the same ``first
come first served'' rule used above. The overall number of
incidences that are lost in this process is again $O(nD)$. The
overall number of polynomials is $O(D^2)$, as can easily be checked.
Note also that the last decomposition step preserves non-singularity
of the points in the special sense defined above; that is, as is
easily verified, a point $p\in Z(f_j)$ with $f_{j+1}(p)\ne 0$,
continues to be a non-singular point of the irreducible component it
is reassigned to.

We now fix one such final polynomial, call it $f_j$, denote its
degree by $D_j$ (which is upper bounded by the original degree $D$),
and denote by $P_j$ and $L_j$ the subsets of the original sets of
points and lines that are assigned to $f_j$, and by $m_j$ and $n_j$
their respective cardinalities. We now may assume that $P_j$
consists exclusively of \emph{non-singular} points of the
\emph{irreducible} variety $Z(f_j)$. The preceding analysis yields
the bound~(\ref{ma:in}) for each $j$
\begin {equation*}
I(P_j,L_j) = O\left( m_j^{1/2}n_j^{1/2}D_j +
m_j^{2/3}n_j^{1/3}s^{1/3} + n_jD_j + m_j\right) .
\end {equation*}
Summing these bounds, upper bounding $D_j$ by $D$, and using
H\"older's inequality for the first two terms, we get the bound~(\ref{ma:in}) for the entire
sets $P$ and $L$. This completes the proof for the case where the
variety containing $P$ and the lines of $L$ are embedded in
$\reals^4$.

\paragraph{Reduction to the four-dimensional case.}
To complete the analysis, we need to consider the
case where $V$ is a three-dimensional variety embedded in
$\reals^d$, for $d>4$. The analysis follows closely the one at the end
of the proof of Theorem~\ref{th:main3d}, in Section~\ref{ssec:pfa}.

Concretely, let $H$ be a generic 4-flat, and denote by $P^*, L^*$, and $V^*$ the
respective projections of $P, L$, and $V$ onto $H$. Since $H$ is
generic, we may assume that all the projected points in $P^*$ are
distinct, and so are all the projected lines in $L^*$. Clearly,
every incidence between a point of $P$ and a line of $L$ corresponds
to an incidence between the projected point and line. Since no
2-flat contains more than $s$ lines of $L$, and $H$ is generic,
repeated applications of Lemma~\ref{ap:th} imply that no 2-flat in
$H$ contains more than $s$ lines of $L^*$.

As in Section~\ref{ssec:pfa}, the set-theoretic projection $V^*$ of
$V$ does not have to be a real algebraic variety, so we use instead
the algebraic projection $\tilde V$ of $V$ that contains $V^*$. That
$\tilde V$ does not contain a hyperplane or quadric follows by a
suitable adaptation of the preceding argument (see~\cite[Lemma
2.1]{SSsocg}). The case of a hyperplane is straightforward
(reasoning as in the preceding section). For quadrics we have:

\medskip

\noindent{\bf Claim.}
Let $X$ be a three-dimensional real algebraic variety in $\reals^d$, 
for $d\ge 5$, such that a generic (real) algebraic projection of $V$ 
on $\reals^4$ is a quadric. Then $X$ is a quadric in $\reals^d$.

\noindent{\bf Proof.} Assume to the contrary that $X$ is not a
quadric in $\reals^d$. This implies that, for a generic $(d-2)$-flat
$h$, the curve $C_h = X\cap h$ is not a quadratic curve. For any
2-flat $g\subset h$, let $C_{h,g}$ denote the projection of $C_h$
onto $g$. This implies that for a generic choice of $g$ and a
$(d-2)$-flat $h$ satisfying $g\subset h \subset \reals^d$, the curve
$C_{h,g}$ is not a quadratic planar curve (that is, a conic section)
in $g$. Next, by taking a suitable rotation of the coordinate
frame, we may assume that $g$ is the $x_1x_2$-flat, and $h$ is the
$x_1x_2\ldots x_{d-2}$-flat. In these coordinates, it is easy to
verify that $C_{h,g}$ can be obtained by first projecting $X$ onto
the $x_1x_2x_{d-1}x_d$-flat, and then cutting it with the
$x_1x_2$-plane. But the projection of $X$ onto the
$x_1x_2x_{d-1}x_d$-flat (which is actually a generic 4-flat) is a
3-quadric by assumption, and then cutting it with any 2-flat
yields a quadratic planar curve, a contradiction that completes the
proof. \proofend

In conclusion, we have $I(P,L) \le I(P^*,L^*)$, where $P^*$ is a set
of $m$ points and $L^*$ is a set of $n$ lines, all contained in the
three-dimensional algebraic variety $\tilde V$, embedded in 4-space,
which is of degree at most $D$ and does not contain any hyperplane
or quadric component, and no 2-flat contains more than $s$ lines of $L^*$.
The preceding analysis thus implies that the bound asserted in the
theorem applies in any dimension $d\ge 4$. \proofend

\section{Discussion}
\label{Se:dis}

{\bf (1)} As already emphasized, most of the analysis in the proof of
Theorem~\ref{th:main3d}(a) is carried out over the complex domain.
The only place where the proofs of (a) and (b) bifurcate is in the final step.
Over the reals we bound $I(P,L_0)$ using the bound of Guth and Katz (which only
holds over the reals, because of the polynomial partitioning that it employs), 
as a ``black box''. Over the complex field, we use a variant of the analysis
to bypass this step, and obtain more or less the same bound, except for 
the additional term $O(D^3)$, which becomes insignificant for $D=O(m^{1/3})$, say.

For three-dimensional varieties, the proof of both parts of
Theorem~\ref{th:main4d} are more or less the same, with the main difference being
the application of the real or complex version of Corollary~\ref{th:degsquared}.
Another difference is in the application of the planar point-line incidence bound---the 
bound is the same in both cases, but the sources (Szemer\'edi-Trotter or T\'oth-Zahl) are different.

The derivation of fairly sharp point-line incidence bounds over the complex domain 
in higher dimensions constitutes, in our opinion, significant progress in this theory.

\medskip

\noindent{\bf (2)}
In view of the lower bound constructions in \cite{GK2,SS4d,SZ}, the new bounds 
in Theorems~\ref{th:main3d} and \ref{th:main4d} do not hold without the assumption
that the points lie on a variety of relatively small degree. We also note that, 
for a three-dimensional variety, we also get rid of the term $m^{1/2}n^{1/2}q^{1/4}$; 
this term may arise only when we consider points on hyperplanes or quadrics, but in 
our case the variety does not contain any such components. Therefore, our theorems 
indicate that these terms may \emph{only} arise if the variety contains such components.

\medskip

\noindent{\bf (3)}
As mentioned in the introduction, Corollary~\ref{cor:main} can be extended to
the case where $V$, which is of constant degree $\Deg$, also contains planes.
Here too, we assume that no plane contains more than $s$ lines of $L$, but
this time it is not necessarily the case that $s\le \Deg$.

Let $\pi_1,\ldots,\pi_k$ denote the planar components of $V$, where $k\le\Deg=O(1)$.
For each $i=1,\ldots,k$, the number of incidences within $\pi_i$, namely, between the
set $P_i$ of points contained in $\pi_i$ and the set $L_i$ of lines fully contained
in $\pi_i$, in both real and complex cases, is
$$
I(P_i,L_i) = O\left( |P_i|^{2/3}|L_i|^{2/3} + |P_i| + |L_i| \right)
= O\left( m^{2/3}s^{2/3} + m + s \right).
$$
Summing these bounds over the $k=O(1)$ planes, we get the same asymptotic bound for
the overall number of the incidences within these planes. Any other incidence
between a point $p$ lying in one of these planes $\pi_i$ and a line $\ell$
not contained in $\pi_i$ can be uniquely identified with the intersection of $\ell$
with $\pi_i$. The overall number of such intersections is at most $nk = O(n)$.
This leads to the following extension of Corollary~\ref{cor:main}.
\begin{corollary}
\label {cor:mainx} Let $P$ be a set of $m$ distinct points and $L$ a
set of $n$ distinct lines in $\reals^d$ or in $\cplx^d$, for any $d\ge 3$, 
and let $s\le n$ be a parameter,
such that all the points and lines lie in a common two-dimensional algebraic surface of
constant degree, and no 2-flat contains more than $s$ lines of $L$. Then
$$
I(P,L) = O\left(m^{2/3}s^{2/3} + m + n\right) ,
$$
where the constant of proportionality depends on the degree of the surface.
\end{corollary}

\medskip

\noindent{\bf (4)} 
As already noted, one of the significant achievements of the
analysis in Theorem~\ref{th:main3d} is that the bound there does not
include the term $O(n\Deg)$. Such a term arises naturally, 
when one considers incidences between points lying in some
irreducible component of $V$ and lines not contained in that component. 
These incidences can be bounded by $n\Deg$, by charging
them, as above, to line-component intersections. 
When $D$ is large, eliminating the term $nD$ can be crucial for the analysis, 
as demonstrated in our earlier work~\cite{SS4d}.

\medskip

\noindent{\bf (5)} Another interesting challenge is to establish a
similar bound for $I(P,L)$, for the case where the points of $P$ lie
on a two-dimensional variety $V$, but the lines need not be
contained in $V$. A trivial extension of the proof adds the term $O(nD)$
to the bound. The challenge is to avoid this term (if possible);
see also Remark (4) above.

\medskip

\noindent{\bf (6)} 
Similar to item (2), Theorem~\ref{th:main4d}(a) can be 
extended to the case where $V$ also contains hyperplane and quadric components,
albeit only for the real case.
Here, as in Theorem~\ref{th:mainfocs}, we add the condition that no
hyperplane or quadric contains more than $q$ lines of $L$.

Let $H_1,\ldots,H_k$ denote the hyperplane and quadric components
of $V$, where $k\le\Deg$. Assign, whenever applicable, each point (resp., line) to the
first $H_i$ that (fully) contains it. As observed above, the number
of ``cross-incidences'' is $O(nD)$. By~\cite[Proposition 3.6]{SS4d},
the total number of incidences within the hyperplanes and quadrics
$H_i$, for $i=1,\ldots, k$, is
$$
O(m^{1/2}n^{1/2}q^{1/4} + m^{2/3}n^{1/3}s^{1/3}+m+n).
$$
This leads to the following extension of Theorem~\ref{th:main4d}(a).
\begin{corollary}
\label {cor:main4dx} Let $P$ be a set of $m$ distinct points and $L$
a set of $n$ distinct lines in $\reals^d$, for any $d\ge 4$, and let $s \le q \le n$
be parameters, such that all the points and lines lie in a common
three-dimensional algebraic surface of degree $D$, and assume that (i) no 3-flat or
3-quadric contains more than $q$ lines of $L$, and (ii) no 2-flat
contains more than $s$ lines of $L$. Then
\begin {equation}
\label{ma:in00} I(P,L) = O\left( m^{1/2}n^{1/2}(D + q^{1/4}) +
m^{2/3}n^{1/3}s^{1/3} + nD + m\right) .
\end {equation}
\end{corollary}

\medskip

\noindent{\bf (7)} An interesting offshoot of Lemma~\ref{claim:4D}
is the following result.

\begin{proposition} \label{2rich}
Let $V$ be a possibly reducible two-dimensional
algebraic surface of degree $D>1$ in $\reals^3$ or in $\cplx^3$,
with no plane or regulus components, and let $L$ be a set of $n$ distinct 
lines fully contained in $V$. Then the number of $2$-rich points
(points incident to at least two lines of $L$) is $O(nD)$. 
\end{proposition}
\noindent{\bf Proof.}
Partition $L$ into the subsets $L_1$ and $L_0$, as in the proof of
Theorem~\ref{th:main3d}. (Recall that $L_0$ is the set of all lines that are
either contained in non-ruled components of $V$, or contained in more than 
one component, or are exceptional lines on ruled components.)
By Lemma~\ref{claim:4D}, each line of $L_1$
is non-conically incident to only $O(D)$ other lines of $L_1$, for a total
of $O(nD)$ $2$-rich points of this sort. Note that we now carry out the analysis
without pruning any point (we do not want to do that), because $V$ 
does not contain any plane or regulus component.

The number of 2-rich points that are exceptional points is at most the
number of irreducible components of $V$, that is, at most $D$, so this 
number is negligible.

The number of lines in $L_0$ is $O(D^2)$. Let $h$ be a plane or a regulus.
The number of lines of $L_0$ contained in $h$ is at most $\deg(V\cap h) \le 2D = O(\sqrt{|L_0|})$
(this holds if we assume, as we may, that $|L_0| = \Theta(D^2)$).
It therefore follows from Guth and Katz~\cite{GK2} that the number of 2-rich 
points involvoing the lines of $L_0$ is $O(|L_0|^{3/2}) = O(|L_0|D)$.

It remains to consider 2-rich points that are intersection points of a line in $L_0$
and a line in $L_1$. By construction, each line $\ell \in L_1$ is contained 
in precisely one (ruled) component $W$ of $V$. If $p\in \ell$ is also incident to 
a line $\ell'\in L_0$ then, again by construction, $\ell'$ is fully
contained in another component $W'$ of $V$, which does not fully contain $\ell$.
Hence $\ell$ intersects $W'$ in at most $\deg(W')$ points (one of which is $p$), 
for a total of at most $\deg(V)=D$ points. Therefore, the number of 2-rich points involving one
line in $L_1$ and another in $L_0$ is at most $n_1 D$.

As we have exhausetd all cases, the assertion follows.
$\Box$

\medskip

\noindent{\bf (8)} Challenging directions for further research are 
(a) to bound the number of incidences between points and lines on 
$(d-1)$-dimensional varieties in $\reals^d$ (or in higher dimensions), for $d\ge 5$,
(b) to bound the number of 2-rich points in a finite set of lines contained in such a variety, and
(c) to bound the number of incidences between points on a variety and $k$-flats 
(under suitable restrictions) in three, four, or higher dimensions.


\section{Algebraic tools and ruled surfaces} \label{sec:ruled} 

In this section we review the preliminary algebraic (and differential)
geometry infrastructure needed for our analysis, and then go on to establish 
the properties of ruled surfaces that we will use. These properties are 
considered folklore in the literature; having failed to find rigorous proofs 
of them (except for several short proofs or proof sketches for some of them), we
provide here such proofs for the sake of completeness. Some of the
notions covered in this section are also discussed in our
study~\cite{SS4d} on point-line incidences in four dimensions.

\paragraph{Singularity.} The notion of
singularities is a major concept, treated in full generality in
algebraic geometry (see, e.g., Kunz~\cite[Theorem VI.1.15]{Kunz} and
Cox et al.~\cite{CLO}). Here we only recall some of their
properties, and only for a few special cases that are relevant to
our analysis.

Let $V$ be a two-dimensional variety in $\reals^3$ or $\cplx^3$ of
degree $D$, given as the zero set $Z(f)$ of some trivariate
polynomial $f$. Assuming $f$ to be square-free, a point $p\in Z(f)$
is \emph{singular} if $\nabla f(p)=0$. For any point $p\in Z(f)$,
let
$$
f(p+x)=f_{\mu}(x)+f_{\mu+1}(x)+\ldots
$$
be the Taylor expansion of $f$ near $p$, where $f_j$ is the $j$-th
order term in the expansion (which is a homogeneous polynomial of
$x$ of degree $j$), and where we assume that there are no terms of
order (i.e., degree) smaller than $\mu$. (The terms $f_j$ also depend on $p$,
which we regard as fixed in the present discussion.) In general, we have 
$f_1(x) = \nabla f(p)\cdot x$, $f_2(x) = \tfrac12 x^TH_f(p)x$, where $H_f$
is the Hessian matrix of $f$, and the higher-order terms are
similarly defined, albeit with more involved expressions.

If $p$ is singular, we have $\mu\ge 2$. In this case, we say that
$p$ is a singular point of $V=Z(f)$ of \emph{multiplicity}
$\mu=\mu_V(p)$. For any point $p\in Z(f)$, we call the hypersurface
$Z(f_{\mu})$ the \emph{tangent cone} of $Z(f)$ at $p$, and denote it
by $C_p Z(f)$. If $\mu=1$, then $p$ is non-singular and the tangent
cone coincides with the (well-defined) tangent plane $T_p Z(f)$ to
$Z(f)$ at $p$. We denote by $V_{sing}$ the locus of singular points
of $V$. This is a subvariety of dimension at most $1$; see, e.g.,
Solymosi and Tao~\cite[Proposition~4.4]{SoTa}. We say that a line
$\ell$ is a \emph{singular line} for $V$ if all of its points are
singular points of $V$.

Similarly, let $\gamma$ be a one-dimensional algebraic curve in
$\reals^2$ or in $\cplx^2$, specified as $Z(f)$, for some bivariate
square-free polynomial $f$. Then $p\in Z(f)$ is \emph{singular} if
$\nabla f(p) = 0$. The \emph{multiplicity} $\mu$ of a point $p\in
\gamma$ is defined as in the three-dimensional case, and we denote
it as $\mu_{\gamma}(p)$; the multiplicity is at least $2$ when $p$
is singular. The singular locus $\gamma_{sing}$ of $\gamma$ is now a
discrete set. Indeed, the fact that $f$ is square-free guarantees
that $f$ has no common factor with any of its first-order
derivatives, and B\'ezout's Theorem (see, e.g., \cite[Theorem
8.7.7]{CLO}) then implies that the common zero set of $f$, $f_x$,
$f_y$, and $f_z$ is a (finite) discrete set.

Still in two dimensions, a line $\ell$, not fully contained in the curve
$\gamma$, can intersect it in at most $D$ points,
\emph{counted with multiplicity}. To define this concept formally,
as in, e.g., Beltrametti~\cite[Section 3.4]{BCGB}, let $\ell$ be a
line and let $p\in\ell\cap\gamma$, such that $\ell$ is not contained
in the tangent cone of $\gamma$ at $p$.
The \emph{intersection multiplicity} of $\gamma$ and $\ell$ at $p$
is the smallest order of a nonzero term of the Taylor expansion of
$f$ at $p$ in the direction of $\ell$. As it happens, the intersection
multiplicity is also equal to $\mu_\gamma(p)$ 
(informally, this is the number of branches of $\gamma$ that $\ell$ crosses at $p$,
counted with multiplicity; see \cite[Section 8.7]{CLO} for a
treatment on the intersection multiplicity in the plane).
The intersection between a line $\ell$ and a curve $\gamma$ (not containing $\ell$)
consists of at most $\deg(\gamma)$ points, counted with their intersection multiplicities.

Assume that $V$ is irreducible. By Guth and Katz~\cite{GK} (see also
Elekes et al.~\cite[Corollary 2]{EKS}), the number of singular lines
contained in $V$ is at most $D(D-1)$.

\paragraph{Flatness.}
We say that a non-singular point $x\in V$ is \emph{flat} if the
second-order Taylor expansion of $f$ at $x$ vanishes on the tangent
plane $T_x V$, or alternatively, if the \emph{second fundamental
form} of $V$ vanishes at $x$ (see, e.g., Pressley~\cite{Pr}).
As argued, e.g., in Elekes et al.~\cite{EKS}, if $x$ is a non-singular point 
of $V$ and there exist three lines incident to $p$ that are fully contained in $V$ 
(this property is captured by calling $p$ a \emph{linearly flat} point) then 
$x$ is a flat point. Following Guth and Katz~\cite{GK}, 
Elekes et al.~\cite[Proposition 6]{EKS} proved that a non-singular point 
$x \in V$ is flat if and only if certain three polynomials, each of degree at 
most $3D-4$, vanish at $p$. A non-singular line $\ell$ is said to be \emph{flat}
if all of its non-singular points are flat. By Guth and Katz~\cite{GK} (see also 
Elekes et al.~\cite[Proposition 7]{EKS}), the number of flat lines fully contained 
in $V$ is at most $D(3D-4)$, unless $V$ is a plane.

As in the proof of Theorem~\ref{th:main4d}, the notions of linear flatness and flatness
can be extended to any higher dimension. For example, for a three-dimensional surface 
$V$ in $\reals^4$ or in $\cplx^4$, which is the zero set of some polynomial $f$ of degree $D$,
a non-singular point $x\in V$ is said to be \emph{linearly flat}, if it is
incident to at least three distinct 2-flats that are fully contained
in $V=Z(f)$ (and thus also in the tangent hyperplane $T_p Z(f)$).
Linearly flat points can then be shown to be \emph{flat}, meaning
that the second fundamental form of $f$ vanishes at them.
This property, at a point $p$, can be expressed by several polynomials of degree $3D-4$
vanishing at $p$ (see~\cite[Section 2.5]{SS4d}).
As in the three-dimensional case, 
the second fundamental form vanishes identically on $Z(f)$ if and only if
$Z(f)$ is not a hyperplane. This property holds in any dimension;
see, e.g.~\cite[Exercise 3.2.12.2]{IL}).
As in three dimensions, we call a line fully contained in $V$ \emph{flat} 
if all its non-singular points are flat.

\paragraph{Ruled surfaces.} 
For a modern approach to ruled surfaces, there are many references; see, e.g.,
Hartshorne~\cite[Section V.2]{Hart83}, or Beauville~\cite[Chapter III]{Beau}; 
see also Salmon~\cite{salmon} and Edge~\cite{Edge} for earlier treatments of 
ruled surfaces. Three relevant very recent additions are the survey~\cite{Gut:surv}
and book~\cite{Gut:book} of Guth, as well as a survey in Koll\'ar~\cite{Kollar},
where this topic is addressed in detail.

We say that a real (resp., complex) surface $V$ is \emph{ruled by real} 
(resp., \emph{complex}) \emph{lines} if every point $p$ in a Zariski-open\footnote{%
  The Zariski closure of a set $Y$ is the smallest (by containment) algebraic variety
  $V$ that contains $Y$. $Y$ is Zariski closed if it is equal to its
  closure (and is therefore a variety), and is (relatively) Zariski open
  if its complement (within a given variety) is Zariski closed.
  See Cox et al.~\cite[Section 4.2]{CLO} for further details.}
dense subset of $V$ is incident to a real (complex) line that is
fully contained in $V$. This definition has been used in several
recent works, see, e.g.,~\cite{GK2,Kollar}; it is a slightly weaker
condition than the classical condition where it is required that
\emph{every} point of $V$ be incident to a line contained in $V$
(e.g., as in~\cite{salmon}). Nevertheless, similarly to the proof of
Lemma 3.4 in Guth and Katz~\cite{GK2}, a limit argument implies that
the two definitions are in fact equivalent. We give, in
Lemma~\ref{le:rs} below, a short algebraic proof of this fact, for
the sake of completeness.

\paragraph{Flecnodes in three dimensions and the Cayley-Salmon-Monge Theorem.}
We first recall the classical theorem of Cayley and Salmon, also due to Monge. Consider
a polynomial $f \in \cplx[x,y,z]$ of degree $D\ge 3$. A
\emph{flecnode} of $f$ is a point $p\in Z(f)$ for which there exists
a line that is incident to $p$ and \emph{osculates} to $Z(f)$ at $p$
to order three. That is, if the direction of the line is $v$ then
$f(p) = 0$, and $\nabla_v f(p) = \nabla_v^2 f(p) = \nabla_v^3 f(p) = 0$, 
where $\nabla_v f, \nabla^2_v f, \nabla^3_v f$ are, respectively, the 
first, second, and third-order derivatives of $f$ in the direction $v$ 
(compare with the definition of singular points, as reviewed earlier, 
for the explicit forms of $\nabla_v f$ and $\nabla^2_v f$). 
The \emph{flecnode polynomial} of $f$, denoted $\fl{f}$, is the polynomial 
obtained by eliminating $v$ from these three homogeneous equations 
(where $p$ is regarded as a fixed parameter). As
shown in Salmon~\cite[Chapter XVII, Section III]{salmon}, the degree
of $\fl{f}$ is at most $11D-24$. By construction, the flecnode
polynomial of $f$ vanishes on all the flecnodes of $f$, and in
particular on all the lines fully contained in $Z(f)$.

\noindent{\bf Theorem 2.1.} (Cayley and Salmon~\cite{salmon}, Monge~\cite{Mon})
{\it Let $f \in \cplx[x,y,z]$ be a polynomial of degree
$\Deg\ge 3$. Then $Z(f)$ is ruled by (complex) lines if and only if
$Z(f) \subseteq Z(\fl{f})$.}

(Note that the correct formulation of Theorem~\ref{th:flec2a} is
over $\cplx$; earlier applications, over $\reals$, as the one in
Guth and Katz~\cite{GK2}, require some additional arguments to
establish their validity; see Katz~\cite{Katz} for a discussion of
this issue.)

\begin{lemma} \label{le:rs} 
Let $f \in \cplx[x,y,z]$ be an irreducible polynomial such that there exists 
a nonempty Zariski open dense set in $Z(f)$ so that each point in the set is 
incident to a line that is fully contained in $Z(f)$. Then $\fl{f}$ vanishes 
identically on $Z(f)$, and $Z(f)$ is ruled by lines.
\end {lemma}
\noindent{\bf Proof.} 
Let $U\subset Z(f)$ be the set assumed in the lemma. By assumption and definition, 
$\fl{f}$ vanishes on $U$, so $U$, and its Zariski closure, are contained in 
$Z(f,\fl{f})$. Since $U$ is open, it must be two-dimensional. Indeed, otherwise 
its complement would be a (nonempty) two-dimensional subvariety of
$Z(f)$ (a Zariski closed set is a variety). In this case, the
complement must be equal to $Z(f)$, since $f$ is irreducible, which
is impossible since $U$ is nonempty. Hence $Z(f,\fl{f})$ is also
two-dimensional, and thus, by the same argument just used, must be
equal to $Z(f)$. Theorem~\ref{th:flec2a} then implies that $Z(f)$ is
ruled by (complex) lines, as claimed. \proofend

The notions of flecnodes and of the flecnode polynomial can be extended to
four dimensions, as done in~\cite{SS4d}. Informally, the 
\emph{four-dimensional flecnode polynomial} $\fl{f}^4$ of $f$ is
defined analogously to the three-dimensional variant $\fl{f}$, 
and captures the property that a point on $Z(f)$ is incident to 
a line that osculates to $Z(f)$ up to the \emph{fourth} order.
It is obtained by eliminating the direction $v$ of the osculating
line from the four homogeneous equations given by the vanishing of
the first four terms of the Taylor expansion of $f(p+tv)$ near $p$.
Clearly, $\fl{f}^4$ vanishes identically on every line that is fully contained
in $Z(f)$. As in the three-dimensional case, its degree can be shown to be $O(D)$.

Landsberg~\cite{Land} derives an analog of Theorem~\ref{th:flec2a} that holds 
for three-dimensional surfaces (see~\cite[Theorem 2.11]{SS4d}). Specifically, 
Landsberg's theorem asserts that if $\fl{f}^4$ vanishes identically on $Z(f)$, 
then $Z(f)$ is ruled by (possibly complex) lines. We will discuss this
in more detail in Section~\ref{sec:pfmain}. These theorems, in three and 
four dimensions, play an important role in the proofs of the main theorems.

\paragraph{Theorem of the fibers and related tools.}
The main technical tool for the analysis is the following so-called
\emph{Theorem of the Fibers}. Both Theorem~\ref{th:harr} and
Theorem~\ref{th:harrdim} hold (only) for the complex field $\cplx$.

\begin {theorem} [Harris~\protect{\cite[Corollary 11.13]{Har}}]
\label {th:harr} Let $X$ be a projective variety and $\pi: X \to
\P^d$ be a homogeneous polynomial map (i.e., the coordinate
functions $x_0\circ \pi,\ldots,x_d \circ \pi$ are homogeneous
polynomials); let $Y=\pi(X)$ denote the image of $X$. For any $p\in
Y$, let $\lambda(p)=\dim(\pi^{-1}(\{p\}))$. Then $\lambda(p)$ is an
upper semi-continuous function of $p$ in the Zariski topology on
$Y$; that is, for any $m$, the locus of points $p\in Y$ such that
$\lambda(p)\ge m$ is Zariski closed in $Y$. Moreover, if $X_0
\subset X$ is any irreducible component, $Y_0=\pi(X_0)$ its image,
and $\lambda_0$ the minimum value of $\lambda(p)$ on $Y_0$, then
$$
\dim(X_0)=\dim(Y_0)+\lambda_0 .
$$
\end {theorem}

We also need the following theorem and lemma from Harris~\cite{Har}.
\begin {theorem} [Harris~\protect{\cite[Proposition 7.16]{Har}}]
\label{th:harrdim} Let $f: X \to Y$ be the map induced by the
standard projection map $\pi: \P^d \to \P^{r}$ (which retains $r$ of
the coordinates and discards the rest), where $r<d$, where $X \subset \P^d$ 
and $Y\subset \P^{r}$ are projective varieties, $X$ is irreducible,
and $Y$ is the image of $X$ (which is also irreducible). 
Then the general fiber\footnote{%
  The meaning of this statement is that the assertion holds for
  the fiber at any point outside some lower-dimensional
  exceptional subvariety.}
of the map $f$ is finite if and only if $\dim(X)=\dim(Y)$. In this
case, the number of points in a general fiber of $f$ is constant.
\end {theorem}

In particular, when $Y$ is two-dimensional (and $d>r\ge 2$ are arbitrary), there exist an integer
$c_f$ and an algebraic curve $\mathcal C_f \subset Y$, such that for
any $y \in Y\setminus \mathcal C_f$, we have $|f^{-1}(y)|=c_f$. With
the notations of Theorem~\ref{th:harrdim}, the set of points $y\in
Y$, such that the fiber of $f$ over $y$ is not equal to $c_f$ is a
Zariski closed proper subvariety of $Y$. For more details, we refer
the reader to Shafarevich~\cite[Theorem II.6.4]{Shaf}, and to
Hartshorne~\cite[Exercise II.3.7]{Hart83}.

\begin {lemma} [Harris~\protect{\cite[Theorem 11.14]{Har}}]
\label {le:harr} Let $\pi: X\to Y$ be a polynomial map between two
projective varieties $X$, $Y$, with $Y=f(X)$ irreducible. Suppose
that all the fibers $\pi^{-1}(\{p\})$ of $\pi$, for $p\in Y$, are
irreducible and of the same dimension. Then $X$ is also irreducible.
\end {lemma}

\paragraph{Reguli.} 
We rederive here the following (folklore) characterization of \emph{doubly ruled} surfaces 
in $\reals^3$ or $\cplx^3$, namely, irreducible algebraic surfaces, each of whose points is 
incident to at least two distinct lines that are fully contained in the surface. 
Recall that a \emph{regulus} is the surface spanned by all lines that meet three pairwise 
skew lines in $3$-space.\footnote{%
  Technically, in some definitions (cf., e.g., Edge~\cite[Section I.22]{Edge})
  a regulus is a one-dimensional family of generator lines of the actual
  surface, i.e., a curve in the Pl\"ucker or Grassmannian space of lines,
  but we use here the alternative notion of the surface spanned by these lines.}
For an elementary proof that a doubly ruled surface must be a
regulus, we refer the reader to Fuchs and Tabachnikov~\cite[Theorem 16.4]{FT}.
Their proof however is analytic and works only over the reals.


\smallskip

\noindent{\bf Lemma 2.2.}
{\it Let $V$ be an irreducible ruled surface in $\reals^3$ or in $\cplx^3$ 
which is not a plane, and let $\C\subset V$ be an algebraic curve, such that 
every non-singular point $p\in V\setminus \C$ is incident to exactly two lines 
that are fully contained in $V$. Then $V$ is a regulus.}\footnote{%
  Over $\reals$, a regulus is either a hyperbolic paraboloid or a one-sheeted hyperboloid.
  Over $\cplx$, balls (equivalent to hyperboloids) and paraboloids (equivalent to hyperbolic
  paraboloids) are also reguli, and are indeed doubly ruled by complex lines. Not all
  quadrics are reguli, though: for example, the cylinder $y=x^2$ is not a regulus.}

\begin{figure}[htb]
\begin{center}
\input{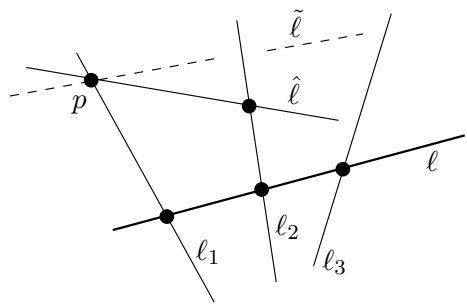}
\caption{The structure of $U_\ell$ in the proof of
Lemma~\ref{doubly}.} \label{freg}
\end{center}
\end{figure}

\noindent{\bf Proof.} 
As mentioned above (see also \cite{GK}), the
number of singular lines in $V$ is finite (it is smaller than
$\deg(V)^2$). For any non-singular line $\ell$, fully contained in
$V$, but not in $\C$, the union of lines $U_{\ell}$ intersecting
$\ell$ and fully contained in $V$ is a subvariety of $V$ (see Sharir
and Solomon~\cite[Lemma 8]{SS3d}
for the easy proof). By assumption, each non-singular point in $\ell \setminus \C$
is incident to another line (other than $\ell$) fully contained in
$V$, and thus $U_\ell$ is the union of infinitely many lines, and is
therefore two-dimensional. Since $V$ is irreducible, it follows that
$U_{\ell}=V$. Next, pick any triple of non-singular and
non-concurrent lines $\ell_1,\ell_2,\ell_3$ that are contained in
$V$ and intersect $\ell$ at distinct non-singular points of $\ell
\setminus \C$. There has to exist such a triple, for otherwise we
would have an infinite family of concurrent (or parallel) lines
incident to $\ell$ and contained in $V$ (where the point of
concurrency lies outside $\ell$), and the plane that they span would
then have to be contained in (the irreducible) $V$, contrary to assumption. 
See Figure~\ref{freg} for an illustration. The argument given for $\ell$
applies equally well to $\ell_1$, $\ell_2$, and $\ell_3$ (by construction, 
neither of them is fully contained in $\C$), and
implies that $U_{\ell_1}=U_{\ell_2}=U_{\ell_3}=V$. 

Assume that there exists some line $\tilde \ell\subset V$ intersecting $\ell_1$ at
some non-singular point $p\in\ell_1\setminus\C$, and that 
$\tilde \ell \cap \ell_2= \emptyset$. (We treat lines here as projective varieties, 
so this assumption means that $\tilde \ell$ and $\ell_2$ are skew to one
another; parallel lines are considered to be intersecting.) Since
$p\in \ell_1 \subset V=U_{\ell_2}$, there exists some line $\hat
\ell$ intersecting $\ell_2$, such that $\hat\ell \cap \ell_1
=\{p\}$. Hence there exist three distinct lines, namely
$\ell_1,\tilde \ell$ and $\hat \ell$, that are incident to $p$ and
fully contained in $V$. Since $p$ is non-singular, it must be a flat point 
(as mentioned above; see~\cite{EKS}). 
Repeating this argument for $3\deg(V)$ non-singular points $p
\in \ell_1$, it follows that $\ell_1$ contains at least $3\deg(V)$
flat points, and is therefore, by the properties of flat points
noted earlier, a flat line. As is easily checked, $\ell_1$ can be
taken to be an arbitrary non-singular line among those incident to
$\ell$, so it follows that every non-singular point on $V$ is flat,
and therefore, as shown in \cite{EKS,GK}, $V$ is a plane, contrary
to assumption.

Therefore, every non-singular line that intersects $\ell_1$ at a
non-singular point also intersects $\ell_2$, and, similarly, it also
intersects $\ell_3$. This implies that the intersection of $V$ and
the surface $R$ generated by the lines intersecting $\ell_1,\ell_2$,
and $\ell_3$ is two-dimensional, and is therefore equal to $V$,
since $V$ is irreducible. Since $\ell_1,\ell_2$ and $\ell_3$ are
pairwise skew, $R=V$ is a regulus, as asserted.  \proofend

\paragraph{Real vs.~complex.}
Expanding on the comment at the end of Section~\ref{Se:pre},
let us elaborate about the field over which the variety $V$ is defined.
Most of the basic algebraic geometry tools have been developed over
the complex field $\cplx$, and some care has to be exercised when
applying them over the reals. A major part of the theory developed
in this section is of this nature. For example, both
Theorems~\ref{th:harr} and~\ref{th:harrdim} hold only over the
complex field. As another important example, one of the main tools
at our disposal is the Cayley--Salmon--Monge theorem (Theorem~\ref{th:flec2a}), 
whose original formulation also applies only over $\cplx$. Expanding on
a previously made comment, we note that even when $V$ is a variety
defined as the zero set of a real polynomial $f$, the vanishing of
the flecnode polynomial $\fl{f}$ only guarantees that the set of
complex points of $V$ is ruled by complex lines.

(A very simple example that illustrates this issue is the unit
sphere $\sigma$, given by $x^2+y^2+z^2=1$, which is certainly not
ruled by real lines, but the flecnode polynomial
of $f(x,y,z)=x^2+y^2+z^2-1$ vanishes on $\sigma$ (since the equation
$\nabla_v^3 f(p) = 0$ that participates in its construction is identically 
zero for any quadratic polynomial $f$). This
is the condition in the Cayley--Salmon--Monge theorem that guarantees that
$\sigma$ is ruled by (complex) lines, and indeed it does, as is
easily checked; in fact, for the same reason, every quadric is ruled
by complex lines.)

This issue has not been directly addressed in Guth and
Katz~\cite{GK2}, although their theory can be adjusted to hold for
the real case too, as noted later in Katz~\cite{Katz}.

This is just one example of many similar issues that one must watch
out for.  It is a fairly standard practice in algebraic geometry
that handles a real algebraic variety $V$, defined by real
polynomials, by considering its complex counterpart $V_\cplx$,
namely the set of complex points at which the polynomials defining
$V$ vanish. The rich toolbox that complex algebraic geometry has
developed allows one to derive various properties of $V_\cplx$, but
some care might be needed when transporting these properties back to
the real variety $V$, as the preceding note concerning the
Cayley--Salmon--Monge theorem illustrates. 
Fortunately, though, passing to the complex domain (and sometimes also to the
projective setting) does not pose any difficulties for deriving upper bounds in incidence 
problems---every real incidence will be preserved, and at worst we will be counting 
additional incidences, on the non-real portion of the extended varieties. 
With this understanding, and with the appropriate caution, we will move freely between
the real and complex domains, as convenient.

We note that most of the results developed in Section~\ref{ssec:pfa} of 
this paper also apply over $\cplx$, except for one crucial step
(where we resort to the application of the result of Guth and Katz~\cite{GK2}, 
which holds only over the reals), due to which 
we do not know how to extend Theorem~\ref{th:main3d} to the complex domain.
Nevertheless, we can derive a weaker variant of it for the complex
case---see a remark to that effect in Section~\ref{Se:dis}.

\paragraph{Lines on a variety.} 
In preparation for the key technical Theorem~\ref{singly}, given below, 
we make the following comments. Lines in three dimensions are parameterized by their
\emph{Pl\"ucker coordinates}, as follows (see, e.g.,
Griffiths and Harris~\cite[Section 1.5]{GrHa}). For two points $x,y
\in \P^3$, given in projective coordinates as $x=(x_0,x_1,x_2,x_3)$
and $y=(y_0,y_1,y_2,y_3)$, let $\ell_{x,y}$ denote the (unique) line
in $\P^3$ incident to both $x$ and $y$. The Pl\"ucker coordinates of
$\ell_{x,y}$ are given in projective coordinates in $\P^5$ as
$(\pi_{0,1},\pi_{0,2},\pi_{0,3},\pi_{2,3},\pi_{3,1},\pi_{1,2})$,
where $\pi_{i,j}=x_iy_j-x_jy_i$. Under this parameterization, the set
of lines in $\P^3$ corresponds bijectively to the set of points in
$\P^5$ lying on the \emph{Klein quadric} given by the quadratic
equation $\pi_{0,1}\pi_{2,3}+\pi_{0,2}\pi_{3,1}+\pi_{0,3}\pi_{1,2} = 0$ 
(which is indeed always satisfied by the Pl\"ucker coordinates of a line). 

Given a surface $V$ in $\P^3$, the set of lines fully
contained in $V$, represented by their Pl\"ucker coordinates in
$\P^5$, is a subvariety of the Klein quadric, which is denoted by $F(V)$, 
and is called the \emph{Fano variety} of $V$; see
Harris~\cite[Lecture 6, page 63]{Har} for details, and \cite[Example 6.19]{Har} 
for an illustration, and for a proof that $F(V)$ is indeed a variety. 
The Pl\"ucker coordinates are continuous, in the sense that if one takes 
two points $\ell$, $\ell'$ on the Klein quadric that are near each other, 
the lines in $\P^3$ that they correspond to are also near to one another, 
in an obvious sense whose precise details are omitted here.

We note again that our analysis is carried out in the complex
projective setting, which makes it simpler, and facilitates the
application of numerous tools from algebraic geometry that are
developed in this setting. The passage from the complex projective
setup back to the real affine one is straightforward---the former is
a generalization of the latter.

Given a plane $\pi$ by a homogeneous equation $A_0 x_0 + A_1 x_1 +
A_2 x_2 + A_3 x_3 = 0$, and a line $\ell$ not fully contained in
$\pi$, given in Pl\"ucker coordinates as
$(\pi_{0,1},\pi_{0,2},\pi_{0,3},\pi_{2,3},\pi_{3,1},\pi_{1,2})$,
their point of intersection is given in homogeneous coordinates by
$(A\cdot m, A \times m - A_0 d)$,
where $d=(\pi_{0,1},\pi_{0,2},\pi_{0,3})$,
$m=(\pi_{2,3},\pi_{3,1},\pi_{1,2})$, and where $\cdot$ stands for
the scalar product, and $\times$ for the vector product; see, e.g.,~\cite[p.~29]{schaum}.
This, together with the continuity argument stated above, implies that, if
the Fano variety $F(V)$ is one-dimensional, and $\ell$ is a line
represented by a non-singular point of $F(V)$, then the cross section of
the union of the lines that lie near $\ell$ in $F(V)$ with a generic
plane $\pi$ is a simple arc. When $\ell$ is a singular point of
$F(V)$, then the cross section of the union of the lines that lie
near $\ell$ in $F(V)$ with a generic plane $\pi$ is a union of
simple arcs meeting at $\ell \cap \pi$ where some of these arcs
might appear with multiplicity; the number of these arcs is determined 
by the multiplicity of the singularity of $\ell$.

\paragraph{Singly ruled surfaces.}
Ruled surfaces that are neither planes nor reguli are called
\emph{singly ruled} surfaces (a terminology justified by
Theorem~\ref{singly}, given below). A line $\ell$, fully contained
in an irreducible singly ruled surface $V$, such that every point of
$\ell$ is ``doubly ruled'', i.e., every point on $\ell$ is incident
to another line fully contained in $V$, is called an
\emph{exceptional} line\footnote{In Guth and Katz~\cite{GK2}, a line
$\ell$ fully contained in an irreducible singly ruled surface $V$,
is called \emph{exceptional} if it contains infinitely many ``doubly
ruled'' points, each incident to another line fully contained in
$V$. Our definition appears to be stricter, but, as the proof below
will reveal, the two notions are equivalent.} of $V$. A point $p_V
\in V$ that is incident to infinitely many lines fully contained in
$V$ is called an \emph{exceptional} point of $V$.

The following result is
another folklore result in the theory of ruled surfaces, used in
many studies (such as Guth and Katz~\cite{GK2}).
It justifies the terminology ``singly-ruled surface'', by showing that
the surface is generated by a one-dimensional family of lines, and that each point
on the surface, with the possible exception of points lying on some curve,
is incident to exactly one generator. It also shows that there are only finitely
many exceptional lines; the property that their number is at most two (see~\cite{GK2})
is presented later. We give a detailed and rigorous proof, to make our presentation
as self-contained as possible; we are not aware of any similarly detailed argument
in the literature.

\smallskip

\noindent{\bf Theorem 2.3.}
{\it (a) Let $V$ be an irreducible ruled two-dimensional
surface of degree $\Deg>1$ in $\reals^3$ (or in $\cplx^3$), which is
not a regulus. Then, except for finitely many exceptional lines, the
lines that are fully contained in $V$ are parameterized by an
irreducible algebraic curve $\Sigma_0$ in the Pl\"ucker space
$\P^5$, and thus yield a 1-parameter family of generator lines
$\ell(t)$, for $t\in \Sigma_0$, that depend continuously on the real
or complex parameter $t$. Moreover, if $t_1 \ne t_2$, and $\ell(t_1)
\ne \ell(t_2)$, then there exist sufficiently small and disjoint
neighborhoods $\Delta_1$ of $t_1$ and $\Delta_2$ of $t_2$, such that
all the lines $\ell(t)$, for $t\in \Delta_1\cup \Delta_2$, are
distinct.}

\smallskip

\noindent {\it (b) There exists a one-dimensional curve $\C\subset V$,
such that any point $p$ in $V\setminus\C$ is incident to exactly one
generator line of $V$.}

\noindent{\bf Remark.} For a detailed description of the algebraic
representation of $V$ by generators, as in part (a) of the theorem,
see Edge~\cite[Section II]{Edge}.

\medskip

\noindent{\bf Proof.} Assume first that we are working over $\cplx$.
Consider the Fano variety $F(V)$ of $V$, as defined above. We claim
that all the irreducible components of $F(V)$ are at most
one-dimensional. Informally, if any component $\Sigma_0$ of $F(V)$
were two-dimensional, then the set $\{ (p,\ell) \in V\times F(V)
\mid p\in\ell\}$ would be three-dimensional, so, ``on average'', the
set of lines of $F(V)$ incident to a point $p\in V$ would be
one-dimensional, implying that most points of $V$ are incident to
infinitely many lines that are fully contained in $V$, which can
happen only when $V$ is a plane (or a non-planar cone, which cannot
arise with a non-singular point $p$ as an apex), contrary to
assumption.

To make this argument formal, consider the set (already mentioned
above)
$$
W:=\{(p,\ell)\mid p \in \ell, \ell \in F(V)\}\subset V \times F(V) ,
$$
and the two projections
$$
\Psi_1 : W \to V, \quad \Psi_2: W \to F(V)
$$
to the first and second factors of the product $V\times F(V)$,
respectively.

$W$ can formally be defined as the zero set of suitable homogeneous
polynomials; briefly, with the Pl\"ucker parameterization of lines
in $\P^3$, and putting the point $p$ into homogeneous coordinates, 
the condition $p\in \ell$ can be expressed as the vanishing of two suitable
homogeneous polynomials, and the other defining polynomials of $W$ are
those that define the projective variety $F(V)$. Therefore, $W$ is a
projective variety.

Consider an irreducible component $\Sigma_0$ of $F(V)$ (which is
also a projective variety); put
$$
W_0:=\Psi_2^{-1}(\Sigma_0)=\{(p,\ell) \in W \mid \ell \in
\Sigma_{0}\}.
$$
Since $W$ and $\Sigma_0$ are projective varieties, so is $W_0$. As
is easily verified, $\Psi_2(W_0)=\Sigma_0$ (that is, $\Psi_2$ is
surjective). We claim that $W_0$ is irreducible. Indeed, for any
$\ell \in \Sigma_0$, the fiber of the map $\Psi_2|_{W_0}: W_0 \to
\Sigma_0$ over $\ell$ is $\{(p,\ell)\mid p \in \ell\}$ which is
(isomorphic to) a line, and is therefore irreducible of dimension
one. As $\Sigma_0$ is irreducible, Lemma~\ref{le:harr} implies that
$W_0$ is also irreducible, as claimed.

For a point $p\in \Psi_1(W_0)$, consider the set
$\Sigma_{0,p}=\Psi_1|_{W_0}^{-1}(\{p\})$, put
$\lambda(p)=\dim(\Sigma_{0,p})$, and let $\lambda_0 := \min_{p \in
\Psi_1(W_0)} \lambda(p)$. By the Theorem of the Fibers
(Theorem~\ref{th:harr}), applied to the map $\Psi_1|_{W_0}: W_0 \to
V$, we have
\begin {equation}
\label {eq:fiber} \dim(W_0)= \dim(\Psi_1(W_0))+\lambda_0.
\end {equation}
We claim that $\lambda_0 = 0$. In fact, $\lambda(p) = 0$ for all
points $p\in V$, except for at most one point. Indeed, if
$\lambda(p)\ge 1$ for some point $p\in V$, then $\Sigma_{0,p}$ is
(at least) one-dimensional, and $V$, being irreducible, is thus a
cone with apex at $p$; since $V$ can have at most one apex, the
claim follows.  Hence $\lambda_0=0$, and therefore
\begin {equation} \label{eq:fibcon}
\dim(W_0) = \dim(\Psi_1(W_0)) \le \dim(V) = 2.
\end {equation}
Next, assume, for a contradiction, that $\dim(\Sigma_0)=2$. For a
point (i.e., a line in $\P^3$) $\ell \in \Psi_2(W_0)$, the set
$\Psi_2|_{W_0}^{-1}(\{\ell\})=\{(p,\ell)\mid p \in \ell\}$ is
one-dimensional (the equality follows from the way $W_0$ is defined).
Conforming to the notations in the Theorem of the Fibers, we have
$\mu(\ell) := \dim\left(\Psi_2|_{W_0}^{-1}(\{\ell\}) \right) = 1$,
and thus $\mu_0:=\min_{\ell \in \Psi_2(W_0)} \mu(\ell)=1$. Also, by
assumption, $\dim(\Psi_2(W_0))=\dim(\Sigma_0)=2$. By the Theorem of
the Fibers, applied this time to $\Psi_2|_{W_0}: W_0 \to \Sigma_0$,
we thus have
\begin {equation}
\label {eq:fiber2} \dim(W_0)= \dim(\Psi_2(W_0))+\mu_0 = 3,
\end {equation}
contradicting Equation~(\ref{eq:fibcon}). Therefore, every
irreducible component of $F(V)$ is at most one-dimensional, as
claimed.

Let $\Sigma_0$ be such an irreducible component, and let
$W_0:=\Psi_2^{-1}(\Sigma_0)$, as above. As argued, for every $p\in
V$, the fiber of $\Psi_1|_{W_0}$ over $p$ is non-empty and finite,
except for at most one point $p$ (the apex of $V$ if $V$ is a cone).
Since $W_0$ is irreducible, Theorem~\ref{th:harrdim} implies that
there exists a Zariski open set $\mathcal O\subseteq V$, such that
for any point $p\in \mathcal O$, the fiber of $\Psi_1|_{W_0}$ over
$p$ has fixed cardinality $c_f$. Put $\C:= V\setminus {\mathcal O}$.
Being the complement of a Zariski open subset of the two-dimensional
irreducible variety $V$, ${\C}$ is (at most) a one-dimensional
variety. If $c_f \ge 2$, then, by Lemma~\ref{doubly}, $V$ is a
regulus. Otherwise, $c_f=1$ ($c_f$ cannot be zero for a ruled
surface), meaning that, for every $p\in V\setminus {\C}$, there is
exactly one line $\ell$, such that $(p,\ell)\in W_0$, i.e.,
$\Sigma_0$ contains exactly one line incident to $p$ and contained
in $V$.

Moreover, we observe that the union of lines of $\Sigma_0$ is the
entire variety $V$. Indeed, by Equations (\ref{eq:fiber}) and
(\ref{eq:fibcon}), we have $\dim(W_0)=\dim(\Psi_1(W_0)) = 2$. That is, 
the variety $\Psi_1(W_0)$, which is the union of the lines of $\Sigma_0$, 
must be the entire variety $V$, because it is two-dimensional and is 
contained in the irreducible variety $V$.

To recap, we have proved that if $\Sigma_0$ is a one-dimensional
component of $F(V)$, then the union of lines that belong to
$\Sigma_0$ covers $V$, and that there exists a one-dimensional
subvariety (a curve) $\C\subset V$ such that, for every $p\in
V\setminus {\C}$, $\Sigma_0$ contains exactly one line incident to
$p$ and contained in $V$.

Since $V$ is a ruled surface, some component of $F(V)$ has to be
one-dimensional, for otherwise we would only have a finite number of
lines fully contained in $V$. We claim that there is exactly one
irreducible component of $F(V)$ which is one-dimensional. Indeed,
assume to the contrary that $\Sigma_0,\Sigma_1$ are two (distinct)
one-dimensional irreducible components of $F(V)$. As we observed,
the union of lines parameterized by $\Sigma_0$ (resp., $\Sigma_1$)
covers $V$. Let $\C_0$, $\C_1\subset V$ denote the respective
excluded curves, so that, for every $p\in V\setminus \C_0$ (resp.,
$p\in V\setminus \C_1$) there exists exactly one line in $\Sigma_0$
(resp., $\Sigma_1$) that is incident to $p$ and contained in $V$.

Next, notice that the intersection $\Sigma_0 \cap \Sigma_1$ is a
subvariety strictly contained in the irreducible one-dimensional
variety $\Sigma_0$ (since $\Sigma_0$ and $\Sigma_1$ are two distinct
irreducible components of $F(V)$), so it must be zero-dimensional,
and thus finite.  Let $\C_{01}$ denote the union of the finitely
many lines in $\Sigma_0 \cap \Sigma_1$, and put $\C:=
\C_0\cup\C_1\cup\C_{01}$. For any point $p\in V\setminus \C$, there
are two (distinct) lines incident to $p$ and fully contained in $V$
(one belongs to $\Sigma_{0,p}$ and the other to $\Sigma_{1,p}$).
Lemma~\ref{doubly} (with $\C$ as defined above) then implies that
$V$ is a regulus, contrary to assumption.

In other words, the unique one-dimensional irreducible component
$\Sigma_0$ of $F(V)$ serves as the desired 1-parameter family of
generators for $V$. (The local parameterization of $\Sigma_0$ can be
obtained, e.g., by using a suitable Pl\"ucker coordinate to represent
its lines.) In addition to $\Sigma_0$, there is a finite number of
zero-dimensional components (i.e., points) of $F(V)$. They
correspond to a finite number of lines, fully contained in $V$, and
not parameterized by $\Sigma_0$. Since the union of the lines in
$\Sigma_0$ covers $V$, any of these additional lines $\ell$ is
\emph{exceptional}, since each point on $\ell$ is also incident to a
generator (different from $\ell$), and is thus ``doubly ruled''.

This establishes part (a) of the theorem, when $V$ is defined over $\cplx$.
We remark that Guth and
Katz~\cite[Corollary 3.6]{GK2} argue that there are at most two such
exceptional lines, so there are at most two zero-dimensional
components of $F(V)$. For the sake of completeness, we sketch a
proof of our own of this fact, in Lemma~\ref{only2exc} below.

If $V$ is defined over $\reals$, we proceed as above, i.e., consider
instead the complex variety $V_\cplx$ corresponding to $V$. As we
have just proven, the unique one-dimensional irreducible component
$\Sigma_0$ of $F(V)$ (regarded as a complex variety) is a (complex)
1-parameter family of generators for the set of complex points of
$V$. Since $V$ is real, the (real) Fano variety of $V$ consists of
the real points of $F(V)$, i.e., it is $F(V)\cap \P^5(\reals)$. As
we have mentioned above, the (complex) $F(V)$ is the union of $\Sigma_0$
with at most two other points. If $\Sigma_0|_\reals := \Sigma_0\cap\P^5(\reals)$
were zero-dimensional, the real $F(V)$ would also be discrete, 
as there is only one one-dimensional component $\Sigma_0$,
so $V$ would fully contain only finitely
many (real) lines, contradicting the assumption that $V$ is ruled by
real lines. Therefore, $\Sigma_0|_\reals$ is a one-dimensional
irreducible component of the real Fano variety of $V$. (It is
irreducible, since otherwise the complex $\Sigma_0$ would be
reducible too, as is easily checked.)

Summarizing, we have shown that there exists exactly one irreducible
one-dimensional component $\Sigma_0$ of $F(V)$, and a corresponding
one-dimensional subvariety ${\C}\subset V$, such that, for each
point $p\in V\setminus {\C}$, $\Sigma_0$ contains exactly one line
that is incident to $p$ (and contained in $V$). In addition to
$\Sigma_0$, $F(V)$ might also contain up to two zero-dimensional
(i.e., singleton) components, whose elements are the exceptional
lines mentioned above. Let $\mathcal D$ denote the union of ${\C}$
and of the at most two exceptional lines; ${\mathcal D}$ is clearly
a one-dimensional subvariety of $V$. Then, for any point $p \in
V\setminus \mathcal D$, there is exactly one line incident to $p$
and fully contained in $V$, as claimed. This establishes part (b),
and thus completes the proof of the theorem. \proofend

\paragraph{Exceptional lines on a singly ruled surface.}
In view of the proofs of Theorem~\ref{singly} and Lemma~\ref{le:rs},
every point on a singly ruled surface $V$ is incident to at least
one generator. Hence an exceptional (non-generator) line is a line
$\ell\subset V$ such that every point on $\ell$ is incident to a
generator (which is different from $\ell$).

\smallskip

\noindent{\bf Lemma 2.4.}
{\it Let $V$ be an irreducible ruled surface in $\reals^3$ or in
$\cplx^3$, which is neither a plane nor a regulus. Then (i) $V$
contains at most two exceptional lines, and (ii) $V$ contains at
most one exceptional point.}

\noindent{\bf Proof.}
(i) We use the property, established in~\cite{SS3d} and already used in the 
proof of Lemma~\ref{doubly}, that for a line $\ell$ fully contained in $V$,
the union $\tau(\ell)$ of the lines that meet $\ell$ and are fully
contained in $V$ is a variety in the complex projective space
$\P^3(\cplx)$. Moreover, if $\ell$ is an exceptional line of $V$,
then it follows by~\cite[Lemma 8]{SS3d} that $\tau(\ell)=V$.
(Indeed, $\tau(\ell)$ must be two-dimensional, since otherwise it
would consist of only finitely many lines. Since $V$ is irreducible,
$\tau(\ell)$ must then be equal to $V$.)

If $V$ contained three exceptional lines, $\ell_1, \ell_2$ and
$\ell_3$, then $V$ would have to be either a plane or a regulus. Indeed,
otherwise, by Theorem~\ref{singly} (whose proof does not depend on
the number of exceptional lines), there would exist a one-dimensional
curve $\C\subset V$ (that includes $\ell_1\cup\ell_2\cup\ell_3$), 
such that every point $p\in V\setminus \C$ is incident to exactly one line 
$\ell_p$ fully contained in $V$. As $p \in V\setminus\C$ and $\sigma(\ell_i)=V$, 
for $i=1,2,3$, it follows that $\ell_p$ intersects $\ell_1,\ell_2$, and $\ell_3$.

If $\ell_1,\ell_2$, and $\ell_3$ are pairwise skew, $p$ belongs to
the regulus $R_{\ell_1,\ell_2,\ell_3}$ of all lines intersecting
$\ell_1,\ell_2,$ and $\ell_3$. We have thus proved that $V\setminus
\mathcal C$ is contained in $R_{\ell_1,\ell_2,\ell_3}$, and as
$R_{\ell_1,\ell_2,\ell_3}$ is irreducible, it follows that $V =
R_{\ell_1,\ell_2,\ell_3}$.

If $\ell_1,\ell_2$, and $\ell_3$ are concurrent but not coplanar
then, arguing similarly, $V$ is a cone with their common
intersection point as an apex. Since a (non-planar) cone has no
exceptional lines, as is easily checked, we may ignore this case.

Finally if any pair among $\ell_1,\ell_2$, $\ell_3$, say $\ell_1$,
$\ell_2$, are parallel then $V$ must be the plane that they span,
contrary to assumption. If $\ell_1$ and $\ell_2$ intersect at a
point $\xi$, disjoint from $\ell_3$, then $V$ is the union of the plane
spanned by $\ell_1$ and $\ell_2$ and the plane spanned
by $\xi$ and $\ell_3$, again a contradiction.

Having exhausted all possible cases, the proof of (i) is complete.

\noindent (ii) 
By Theorem~\ref{singly} and (i), all the lines that are fully contained in $V$, 
except for possibly two such lines, are parameterized by an irreducible algebraic curve
$\Sigma_0$ in the Pl\"ucker space $\P^5$. Let $p$ be an exceptional
point of $V$. The set $\Sigma'$ of lines incident to $p$ is an
algebraic curve contained in the irreducible curve $\Sigma_0$,
implying that $\Sigma'=\Sigma_0$. This clearly implies that there is
at most one exceptional point (and then it does not contain any exceptional line), 
and the proof of (ii) is complete too.
\proofend

\noindent{\bf Remark.} We refer the reader to Guth and
Katz~\cite[Lemma 3.5, Corollary 3.6]{GK2}, for yet another 
(somewhat more compact) proof of this lemma.

\paragraph{Generic projections preserve non-planarity.}
In the analysis in Section~\ref{ssec:pfa}, the goal is to project
$\reals^d$ onto some generic 3-flat so that non-coplanar triples of
lines do not project to coplanar triples. This is easily achieved by
repeated applications of the following technical result, reducing
the dimension one step at a time.

\smallskip

\noindent{\bf Lemma 2.5.}
{\it Let $\ell_1,\ell_2,\ell_3$ be three non-coplanar lines
in $\reals^d$. Then, under a generic projection of $\reals^d$ onto
some hyperplane $H$, the respective images
$\ell^*_1,\ell^*_2,\ell^*_3$ of these lines are still non-coplanar.}

\noindent {\bf Proof.} Assume without loss of generality that the
(generic) hyperplane $H$ onto which we project passes through the
origin of $\reals^d$, and let $w$ denote the unit vector normal to
$H$. The projection $h:\;\reals^d\mapsto H$ is then given by $h(v) =
v-(v\cdot w) w$.

Assume first that two of the three given lines, say $\ell_1,\ell_2$,
are skew (i.e., not coplanar). Let $\tilde \ell_1,\tilde \ell_2$
denote their projection onto $H$. If $\tilde \ell_1, \tilde \ell_2$
are coplanar they are either intersecting or parallel. If they are
intersecting, then there are points $p_1 \in \ell_1, p_2 \in \ell_2$
that project to the same point, i.e., $p_1-p_2$ has the same
direction as $w$. Then $w$ belongs to the set $\{\frac {p_1-p_2}
{\|p_1-p_2\|} \mid p_1 \in \ell_1, p_2 \in \ell_2\}$. Since this is
a two-dimensional set, it will be avoided for a generic choice of
$w$, which is a generic point in $\sph^{d-1}$, a set that is at
least three-dimensional.

If $\tilde \ell_1,\tilde \ell_2$ are parallel, let $v_1,v_2$ denote
the directions of $\ell_1,\ell_2$. Since $v_1-(v_1\cdot w) w$ and
$v_2-(v_2\cdot w) w$ are vectors in the directions of $\tilde
\ell_1,\tilde \ell_2$, and are thus parallel, it follows that $w$
must be a linear combination of $v_1$ and $v_2$. Since $\|w\|=1$,
the resulting set of possible directions is only one-dimensional,
and, again, it will be avoided with a generic choice of $w$.

We may therefore assume that every pair of lines among
$\ell_1,\ell_2,\ell_3$ are coplanar. Since these three lines are not
all coplanar, the only two possibilities are that either they are
all mutually parallel, or all concurrent.

Assume first that they are concurrent, say they all pass through the
origin (even though the origin belongs to $H$, this still involves
no less of generality). Their projections are in the directions
$v_i-(v_i\cdot w) w$, for $i=1,2,3$. If these projections are
coplanar then there exist coefficients $\alpha_1,\alpha_2,\alpha_3$,
not all zero, such that $\sum_i \alpha_i (v_i-(v_i\cdot w) w) = 0$.
That is, putting $u:= \sum_i \alpha_iv_i$, we have $u=(u\cdot w)w$,
so $u$ is parallel to $w$. In this case $w$ belongs to the set
$\left\{\frac{\sum_i \alpha_i v_i} {\|\sum_i \alpha_i v_i\|} \mid
\alpha_1, \alpha_2, \alpha_3 \in \reals \text{ or } \cplx \right\}$.
Again, being a two-dimensional set, it will be avoided by a generic
choice of $w$.

In the remaining case, the lines $\ell_1$, $\ell_2$, $\ell_3$ are
mutually parallel, i.e., they all have the same direction $v$. Put,
for $i=1,2,3$, $\ell_i=\{p_i+tv\}_{t\in\reals}$, and choose $p_i$ so that
$p_i\cdot v = 0$. The plane $\pi_0$ spanned by $p_1,p_2,p_3$ is
projected to the plane $\pi$ spanned by the points $p_i^* = p_i - (p_i\cdot w)w$, 
for $i=1,2,3$ (since $p_1$, $p_2$, $p_3$ are not
collinear, they will not project into collinear points in a generic
projection), and the three lines project into a common plane if and
only if their projections are fully contained in $\pi$, meaning that
the projection $v^* = v-(v\cdot w)w$ is parallel to $\pi$, so it
must be a linear combination of $p_1^*$, $p_2^*$, and $p_3^*$. A
similar argument to those used above shows that a generic choice of
$w$ will avoid the resulting two-dimensional set of forbidden
directions. This completes the proof. \proofend

\end{document}